\documentclass[11pt]{article}
\usepackage[a4paper,margin=25mm]{geometry}
\usepackage[T1]{fontenc}
\usepackage{lmodern,microtype,amsmath,amssymb,amsthm}
\usepackage{pgfplots,flafter}
\pgfplotsset{compat=1.18}
\usepackage[colorlinks=true,linkcolor=blue!42!black,citecolor=blue!42!black,urlcolor=blue!42!black]{hyperref}
\hypersetup{pdftitle={The tower transition in majority bootstrap percolation},pdfauthor={Vyacheslav Koval and Yuval Peres},pdfsubject={Majority bootstrap percolation on high-dimensional tori and boxes}}
\newtheorem{theorem}{Theorem}[section]
\newtheorem{lemma}[theorem]{Lemma}
\newtheorem{corollary}[theorem]{Corollary}
\newtheorem{proposition}[theorem]{Proposition}
\numberwithin{equation}{section}
\newcommand{\T}{\mathcal T}
\newcommand{\ZZ}{\mathbb Z}
\newcommand{\PP}{\mathbb P}
\newcommand{\EE}{\mathbb E}
\newcommand{\Bin}{\operatorname{Bin}}
\title{\vspace{-8mm}The tower transition in majority bootstrap percolation}
\author{Vyacheslav Koval\thanks{Bernoulli Institute, University of Groningen, The Netherlands.
\href{mailto:v.v.koval@rug.nl}{\nolinkurl{v.v.koval@rug.nl}}}
\qquad Yuval Peres\thanks{Beijing Institute of Mathematical Sciences and Applications (BIMSA), China.
\href{mailto:yperes@gmail.com}{\nolinkurl{yperes@gmail.com}}}}
\date{}
\begin{document}
\maketitle
\enlargethispage{\baselineskip}

\begin{abstract}
How large must a high-dimensional torus be for majority bootstrap percolation to fill it at a given density below one half? For a torus of side length $n\ge3$, let $H(n)$ be the smallest $h$ such that a tower of $h$ twos has value at least $n$. If $H(n)/d\to\rho\in[0,1]$, the critical density on the $d$-dimensional torus tends to $(1-\rho)/(2-\rho)$. In particular, the smallest side length $n$ for which the critical density is at most $1/4$ satisfies $H(n)\sim 2d/3$. The same curve holds for ordinary boxes of side at least $d$. For boxes with $n/d\to a>0$ we obtain the limiting density $\max\{1/2,1-2^{-1/a}\}$. We determine the first correction to the tower curve away from $\rho=1$.
\end{abstract}

\section{Introduction}\label{sec:intro}
In $r$-neighbour bootstrap percolation, each vertex of a graph is initially
infected with probability $p$, independently of the others. At each round,
a healthy vertex becomes infected if at least $r$ of its neighbours are
infected. Infected vertices stay infected. Write $A_t$ for the infected set
after $t$ rounds and $A_\infty$ for its union over all $t$. On a finite
graph $G$ the critical density is
\[
 p_c(G,r)=\inf\{p:\PP_p(A_\infty=V(G))\ge 1/2\}.
\]
We say that the graph \emph{fills} when $A_\infty=V(G)$.
Our graphs are the nearest-neighbour torus $\ZZ_n^d=(\ZZ/n\ZZ)^d$,
$n\ge3$, and the box $[n]^d=\{1,\ldots,n\}^d$. On both we take $r=d$,
including at boundary vertices of the box.

A tower of $d$ twos suffices for the critical density to tend to zero,
even on the ordinary box~\cite{CM,BBDCM}. Balogh, Bollob\'as and
Morris~\cite[Problem~1]{BBM} asked how the limiting critical density
depends on the growth of the side length with the dimension.
For density $1/4$, Morris~\cite[\S4.1]{MorrisSurvey} wrote that
\begin{quote}
\textquotedblleft we suspect that the function
$f(d):=\min\{n\in\mathbb N:p_c(\ZZ_n^d,d)\le1/4\}$
grows like a tower of 2s of height $\Theta(d)$.\textquotedblright
\end{quote}

Set $\T_0=1$, $\T_{k+1}=2^{\T_k}$ and
$H(x)=\min\{k\ge0:x\le\T_k\}=\log_2^*x$ for $x\ge1$.
All unmarked logarithms are natural.
\begin{theorem}[The tower transition]\label{thm:phase}
For $n\ge3$, put $\rho=\min\{H(n)/d,1\}$. Then
\begin{equation}\label{eq:curve}
 p_c(\ZZ_n^d,d)=\frac{1-\rho}{2-\rho}+o(1)
\end{equation}
as $d\to\infty$, uniformly in $n$. The same estimate holds for $[n]^d$
whenever $n\ge d$. For any sequence with $n/d\to a\in(0,\infty)$,
\begin{equation}\label{eq:linear-box}
 p_c([n]^d,d)\ \longrightarrow\ \max\{1/2,\,1-2^{-1/a}\}.
\end{equation}
\end{theorem}

For $n=\T_{\lfloor\rho d\rfloor}$ and fixed $0<\rho\le1$, this gives
$p_c\to(1-\rho)/(2-\rho)$ on both graphs. In particular, the critical
density tends to zero when $H(n)\ge d$.

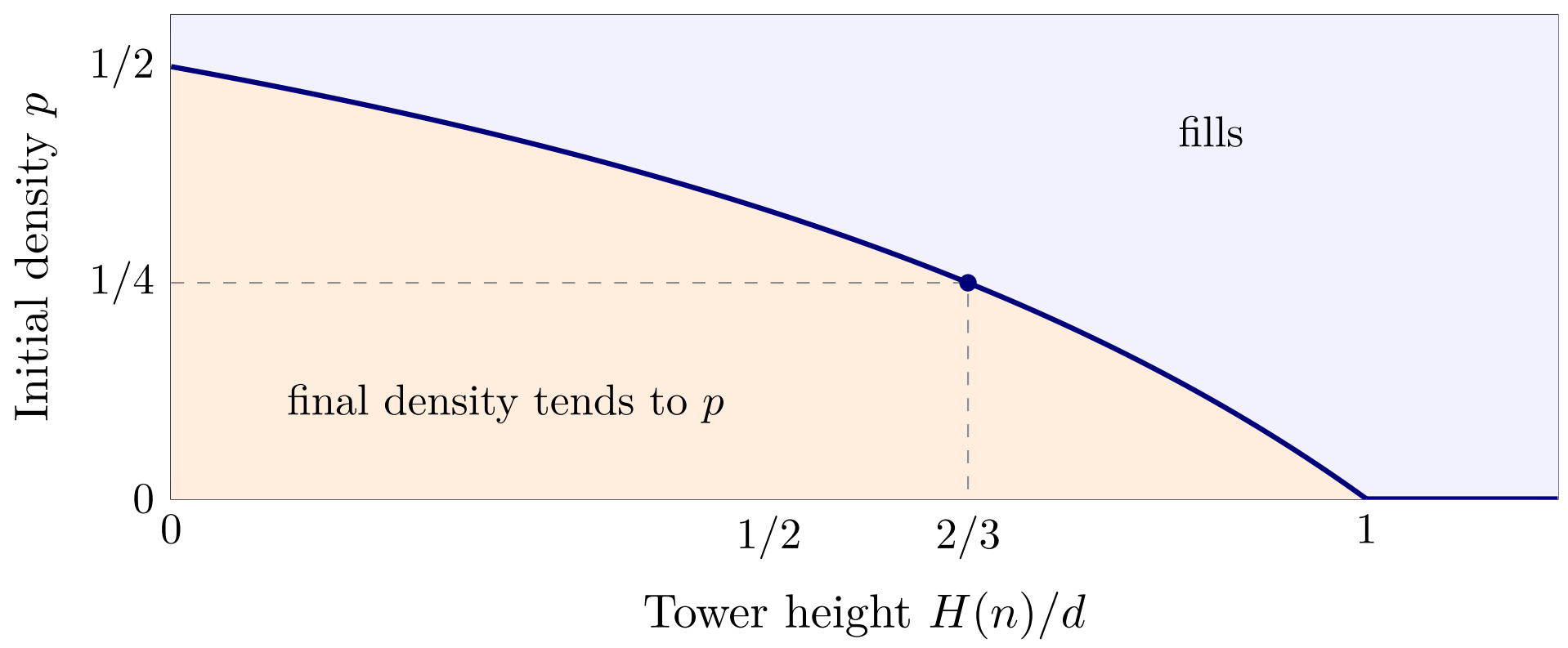
\begin{figure}[htbp]
\centering
\begin{tikzpicture}
\begin{axis}[width=.80\textwidth,height=5.5cm,xmin=0,xmax=1.16,ymin=0,ymax=.56,
 xlabel={Tower height $H(n)/d$},ylabel={Initial density $p$},
 xtick={0,.5,.66666667,1},xticklabels={$0$,$1/2$,$2/3$,$1$},
 ytick={0,.25,.5},yticklabels={$0$,$1/4$,$1/2$},
 tick label style={font=\small},clip=true]
\path[fill=blue!5] (axis cs:0,0) rectangle (axis cs:1.16,.56);
\addplot[draw=none,fill=orange!13,domain=0:1,samples=101] {(1-x)/(2-x)} \closedcycle;
\addplot[draw=blue!48!black,very thick,domain=0:1,samples=101] {(1-x)/(2-x)};
\addplot[draw=blue!48!black,very thick,domain=1:1.16] {0};
\draw[dashed,black!45] (axis cs:0,.25) -| (axis cs:.66666667,0);
\addplot[draw=blue!48!black,fill=blue!48!black,only marks,mark=*,mark size=1.8pt]
 coordinates {(.66666667,.25)};
\node[font=\small] at (axis cs:.87,.425) {fills};
\node[font=\small,align=center] at (axis cs:.28,.11) {final density tends to $p$};
\end{axis}
\end{tikzpicture}
\caption{The tower transition on the torus. For fixed height ratio and fixed
density strictly above the curve the torus fills with high probability.
Strictly below the curve, the fraction of sites that become newly infected
tends to zero. The marked point is Morris's density $1/4$.}
\label{fig:linear-bounds}
\end{figure}

For example, at side $\T_{\lfloor d/2\rfloor}$ the critical density tends
to $1/3$. At side $\T_{\lfloor2d/3\rfloor}$ it tends to $1/4$
(Figure~\ref{fig:linear-bounds}). In the latter example the torus and
the box in fact fill at $p=1/4$ with probability tending to one.
This uses the correction to~\eqref{eq:curve} in Section~\ref{sec:fine}.
The constant $1$ in the linear side condition $n\ge d$ is sharp. For $a<1$ the limit in
\eqref{eq:linear-box} exceeds $1/2$, and for $2\le n=o(d)$ the critical
density tends to one. For instance, when $n\sim d/2$ the limiting
density is $3/4$. At $n=d$ the number of empty corner edge-lines has
a Poisson limit in the critical window (Proposition~\ref{prop:poisson}).

Write $p_\rho=(1-\rho)/(2-\rho)$ for the limiting density.
For $0<p<1/2$ its inverse and the first crossing side length are
\[
 \rho_p=\frac{1-2p}{1-p},\qquad
 f_p(d)=\min\{n\ge3:p_c(\ZZ_n^d,d)\le p\}.
\]
\begin{corollary}[The first crossing]\label{cor:first-crossing}
For each fixed $0<p<1/2$ the minimum exists for all sufficiently large $d$,
and $H(f_p(d))/d\to \rho_p$. In particular, for $0<\varepsilon<2/3$,
\[
 \T_{\lfloor(2/3-\varepsilon)d\rfloor}<f_{1/4}(d)
       \le\T_{\lceil(2/3+\varepsilon)d\rceil}
\]
for all sufficiently large $d$.
\end{corollary}

\begin{theorem}[The final density below the transition]\label{thm:density}
Fix $0<p<1/2$. If $\limsup_{d\to\infty}H(n)/d<\rho_p$, then on the torus
\[
 \frac{|A_\infty\setminus A_0|}{n^d}\ \longrightarrow\ 0
 \quad\hbox{in probability}.
\]
In fact the expected fraction on the left is at most $e^{-cd}$ for
all large $d$, where $c>0$ depends on $p$ and the gap in the hypothesis.
The torus fills with probability at most $e^{-c'd^3}$ for all large $d$,
for a constant $c'>0$ with the same dependence.
The same conclusions hold for the box when $n\ge2$.
\end{theorem}

On either graph the first round accounts for the increase in expected
density up to an exponentially small relative error. On the torus the
leading contribution from later rounds comes from an adjacent pair of
vertices (Proposition~\ref{prop:first-round}).
In the other direction, for fixed $0<p<1/2$ and
$\liminf H(n)/d>\rho_p$, the box and torus both fill with probability
at least $1-e^{-d}$ for all sufficiently large $d$
(Proposition~\ref{prop:leading-up}).

\subsection*{Where the curve comes from}
Suppose $v$ neighbours of a vertex are certain to be infected. Its expected
support, counting initial infections among the other neighbours, is
\[
 v+(2d-v)p.
\]
This reaches $d$ at $v=\rho_p d$. Why should the number of levels in a tower
enter a neighbour count? There are two answers, one for each direction
of the proof.

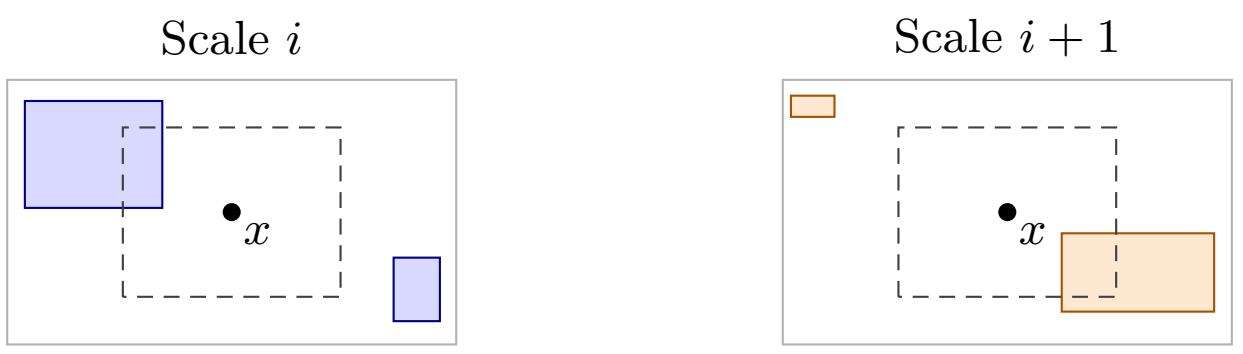
\begin{figure}[htbp]
\centering
\begin{tikzpicture}[x=9mm,y=7mm,font=\footnotesize]
\filldraw[draw=blue!55!black,fill=blue!15] (-1.52,.04) rectangle (-.51,1.05);
\filldraw[draw=blue!55!black,fill=blue!15] (1.19,-1.03) rectangle (1.53,-.43);
\filldraw[draw=orange!65!black,fill=orange!18] (6.10,-.94) rectangle (7.22,-.2);
\filldraw[draw=orange!65!black,fill=orange!18] (4.11,.9) rectangle (4.43,1.1);
\foreach \dx/\layer in {0/i,5.7/{i+1}} {
  \begin{scope}[shift={(\dx,0)}]
  \draw[gray!60] (-1.65,-1.25) rectangle (1.65,1.25);
  \draw[black!75,densely dashed] (-.80,-.80) rectangle (.80,.80);
  \fill (0,0) circle (1.7pt) node[below right,inner sep=2pt] {$x$};
  \node at (0,1.65) {Scale $\layer$};
  \end{scope}
}
\end{tikzpicture}
\caption{The dashed square is the same neighbourhood $B^\infty(x,t)$
at two successive scales. At a given scale the retained rectangles
are more than $2t$ apart, so the square meets at most one. A rectangle
from a different scale may enter the same square.}
\label{fig:one-per-scale}
\end{figure}

For the lower bound we trace an infection outwards through successive
neighbours. A site that gains many neighbours must have many children
with nearly enough initial support. Mark a site if such a configuration
occurs for some choice of the nearby infected rectangles. Merge
rectangles around these rare marks until the small rectangles are well
separated. A large rectangle contains many marks with disjoint input
neighbourhoods, which lets us pass to an exponentially larger scale.
A site outside a retained rectangle has at most one neighbour inside.
Its neighbourhood meets at most one rectangle per scale
(Figure~\ref{fig:one-per-scale}), so each level of the cover supplies at
most one forced neighbour. The local test allows every possible trace
of the rectangles, including the trace of the cover built from the seeds.

For the upper bound, suppose a small cube is already infected. Add one
face. Each site of the face has an infected neighbour behind it, and the
initial state just ahead of it supplies another chance of infection. Finding
a cube from which growth can start costs one exponential in the side length.
Passing to the new face reduces the mean deficit by $1-p$
(Figure~\ref{fig:face}). Repeating these
two steps leads to the same neighbour count. For an ordinary box,
fill its interior first, then its codimension-one faces, and so on down
to the corners. Each free coordinate of the cell being filled has two
neighbours in the original box. Even the last site on a side can use
the initial state of the corner just beyond it.

\subsection*{Earlier work}
Bootstrap percolation was introduced by Chalupa, Leath and Reich~\cite{CLR}.
For fixed dimension and threshold two, Aizenman and Lebowitz~\cite{AL}
used rectangles to find the scale of the critical density. Cerf and
Cirillo~\cite{CC} found the iterated logarithm in dimension three, and
Cerf and Manzo~\cite{CM} treated higher dimensions. The sharp thresholds in
dimensions two and three~\cite{Holroyd,BBM3} were followed by the theorem
of Balogh, Bollob\'as, Duminil-Copin and Morris~\cite{BBDCM}. For fixed
$2\le r\le d$, the critical density on $[n]^d$ is asymptotic to
$(\lambda(d,r)/\log_{(r-1)}n)^{d-r+1}$, where $\log_{(j)}$ denotes $j$
successive logarithms. Their Theorem~28 records the tower upper bound,
and their Problem~2 asks for a box whose critical density stays strictly
between zero and one half. Theorem~\ref{thm:phase} gives a choice for
every limit in this interval. We use the arXiv numbering
for~\cite{BBM,BBDCM,CEGK2}.

The box growth argument has its origins in Schonmann~\cite{Schonmann}
and in the construction of Bou-Rabee, Levine and
Peres~\cite[\S2.11, Lemma~2.45]{Book}, who grow from an infected
corner by adding one face at a time.

The first-order proofs occupy Sections~\ref{sec:lower}
and~\ref{sec:box}. Section~\ref{sec:fine} gives the binomial estimates
that locate the transition more precisely.

\section{The lower bound}\label{sec:lower}
We cover the sites that could be infected late by rectangles. Each scale
can contribute one infected neighbour. A site that can step into a
rectangle from outside has exactly one coordinate outside its interval. On the torus
we use \emph{proper} rectangles, products of cyclic intervals each
omitting at least two sites. Then even in that coordinate only one step
can enter the rectangle. Throughout, side lengths count sites.

Write $B(x,i)$ and $S(x,i)$ for the ball and sphere in graph distance
and $B^\infty(x,i)$ for the ball in toroidal $\ell^\infty$ distance.
We use the following count. For $y\in S(x,i)$,
\begin{equation}\label{eq:backward}
 |N(y)\cap B(x,i)|\le\kappa i,
 \qquad \kappa=2\ \text{on a torus},\quad
 \kappa=1\ \text{in a box or in $\ZZ^d$}.
\end{equation}
Here $N(y)$ is the set of neighbours of $y$. Only coordinates in
which $x$ and $y$ differ can contribute. There are
at most $i$ of them, and each gives at most two choices on a cycle and
one in an interval. This is the backward bound of
Collares, Erde, Geisler and Kang~\cite[Property~P2]{CEGK2}.

\subsection{An outward witness}
Write $d_S(v)=|N(v)\cap S|$ and $q=1-p$.

We use the following form of the Chernoff bound. Suppose $Z$ is a sum
of indicators, each a function of independent bits, and each bit enters
at most $b$ indicators. With $\mu=\EE Z$,
\begin{align}
 \PP(Z\ge u)&\le(e\mu/u)^{u/b} &&(u>\mu>0),\notag\\
 \PP(Z\le(1-\eta)\mu)&\le e^{-\eta^2\mu/(2b)}
                  &&(0<\eta<1).\label{eq:read}
\end{align}
Indeed, successive applications of H\"older's inequality give
$\EE\prod f_i\le\prod(\EE f_i^b)^{1/b}$ for functions with this
overlap~\cite{Finner}. Taking $f_i=e^{\lambda I_i}$ for the indicators
gives $\EE e^{\lambda Z}\le\exp\{\mu(e^{b\lambda}-1)/b\}$.
Optimizing in $\lambda$ gives~\eqref{eq:read}~\cite{GLSS}.

On $\ZZ_n^d$ fix a forced set $V$ and positive integers $h$ and $g$.
Put $G=hg$ and $I=A_0\cup V$. Define
\begin{align}
 E_0(V)&=\{v\notin I:d_I(v)\ge d-G\},\notag\\
 E_i(V)&=\{v:d_{E_{i-1}(V)}(v)>g\}\quad(1\le i\le h),\qquad
 B(V)=I\cup\bigcup_{i=0}^{h-1}E_i(V).\label{eq:marks}
\end{align}
If $E_h(V)\subset B(V)$, then $B(V)$ is closed under the $d$-neighbour
rule. Indeed, a vertex outside it has at most $d-G-1$ neighbours in $I$
and at most $g$ in each of $E_0(V),\ldots,E_{h-1}(V)$. These bounds
sum to $d-1$. In particular $A_\infty\subset B(V)$.

\begin{lemma}[An outward witness]\label{lem:witness}
Suppose $V$ is deterministic and each vertex outside $V$ has at most
$H$ neighbours in $V$, where $H\in\{0,\ldots,2d\}$. Suppose $g>2h$
and put $u=g-2h$. Define
\begin{equation}\label{eq:recurrence}
\begin{aligned}
 w_0&=q\,\PP\{\Bin(2d-H,p)\ge d-G-2h-H\},\\
 w_i&=\min\!\left\{1,\left(\frac{2ed\,w_{i-1}}{u}\right)^{u/[2(i+1)]}\right\}
        \qquad(1\le i\le h).
\end{aligned}
\end{equation}
Then $\PP_p(x\in E_i(V))\le w_i$ for every $x$ and $0\le i\le h$.
\end{lemma}

\begin{figure}[htbp]
\centering
\begin{tikzpicture}[x=1cm,y=1cm,font=\footnotesize,
 v/.style={circle,fill=blue!55!black,inner sep=1.8pt}]
  \node at (0,1.45) {$E_2$};
  \node at (2.35,1.45) {$E_1$};
  \node at (4.9,1.45) {$E_0$};
  \node[v,label=left:$x$] (r) at (0,0) {};
  \node[v] (a) at (2.35,.54) {};
  \node[v] (b) at (2.35,-.54) {};
  \node[v] (c) at (4.9,.9) {};
  \node[v] (e) at (4.9,0) {};
  \node[v] (f) at (4.9,-.9) {};
  \draw[blue!55!black] (r)--(a)--(c) (a)--(e) (r)--(b)--(e) (b)--(f);
  \node[right=4pt] at (e) {a shared site};
  \node at (2.35,-1.32) {$S(x,1)$};
  \node at (4.9,-1.32) {$S(x,2)$};
\end{tikzpicture}
\caption{An outward witness with two generations. Two children of $x$
can use the same mark on $S(x,2)$. A seed one step beyond the leaves can
enter a child test only through a first step on a geodesic to that seed.}
\label{fig:three-rounds}
\end{figure}
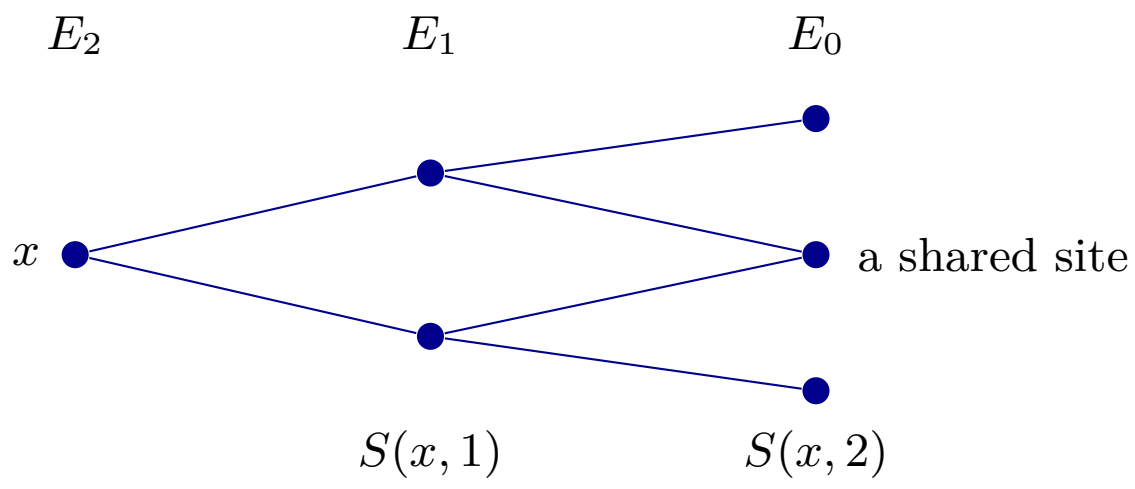

\begin{proof}
Fix $x$ and a depth $j\le h$. A child lies on the next sphere about
$x$. To test $x\in E_j(V)$, work backwards from the definition
(Figure~\ref{fig:three-rounds}). At each step retain only children.
This discards at most $2h$ neighbours by~\eqref{eq:backward}.
A leaf on $S(x,j)$ must be outside $I$ and have at least $d-G-2h$
children in $I$. Replacing up to $H$ of the $2d$ neighbour bits by ones
bounds this probability by $w_0$. The leaf's own independent bit
supplies the factor $q$. Mark a vertex with $i$ generations below it
when at least $u$ children are marked. Membership in $E_j(V)$ implies
a mark at $x$.

Induct within this test, uniformly over its vertices and over $j$.
With $i$ generations left, the child marks use only the bits at the
leaves and one step beyond. A bit at $z$ can enter a test through a
child only if that child is the first step on a geodesic from the
current vertex to $z$. The distance is at most $i+1$.
By~\eqref{eq:backward} with root $z$ there are at most $2(i+1)$
such first steps. The mean number of marked children is at most
$2d w_{i-1}$. The upper tail in~\eqref{eq:read} gives~\eqref{eq:recurrence}.
\end{proof}

A union of at most $H$ proper
rectangles supplies at most $H$ forced neighbours at a vertex outside
the union. There the threshold minus the mean initial support is at least
\begin{equation}\label{eq:deficit}
 D_H=d-[H+(2d-H)p]=(1-2p)d-(1-p)H .
\end{equation}

\begin{lemma}[A linear deficit]\label{lem:three}
Fix $0<p<1/2$ and $\varepsilon>0$. If $H\in\{0,\ldots,d\}$ and
$D_H\ge\varepsilon d$, take $h=2$ and $g=\lfloor\varepsilon d/8\rfloor$.
For large $d$, uniformly in the side and the forced rectangles,
\[
 \PP_p(x\in E_2(V))\le e^{-cd^3},\qquad
 \PP_p(x\in B(\varnothing)\setminus A_0)\le e^{-cd}.
\]
\end{lemma}
\begin{proof}
A leaf needs an excess above its mean of at least
$D_H-2g-4\ge\varepsilon d/2$. Hoeffding's inequality~\cite{Hoeffding}
gives $w_0\le e^{-c_0d}$. Since $u=g-4$ is of order $d$,
two applications of~\eqref{eq:recurrence} give $w_1\le e^{-c_1d^2}$
and $w_2\le e^{-c_2d^3}$.
With empty background the probability of an extra site in $B$ is at
most $w_0+w_1$.
\end{proof}

\subsection{Covering a rare set by separated rectangles}
A random set $Q$ is \emph{$R$-local} if membership of $x$ is determined by
the initial bits in $B^\infty(x,R)$. In one step we cover most marks
by rectangles so far apart that a test window meets at most one.
The remaining marks pass a rarer local test. The radius $t$ will be
three here and six in Section~\ref{sec:fine}.

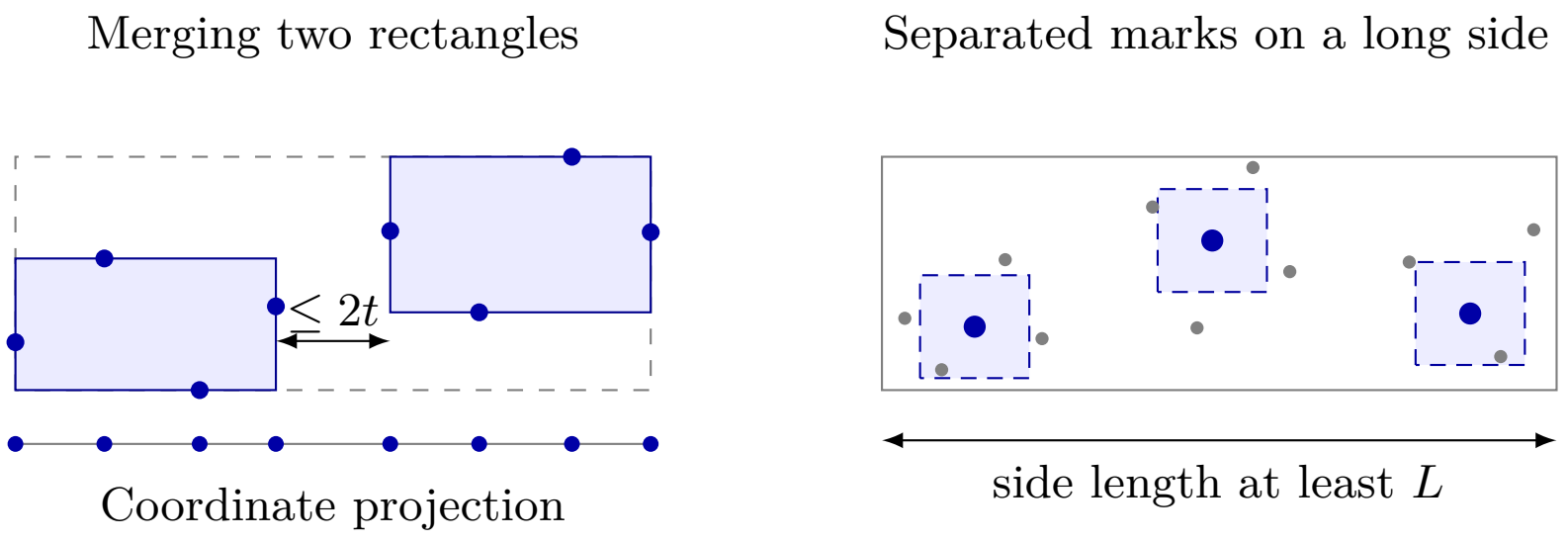
\begin{figure}[htbp]
\centering
\begin{tikzpicture}[x=.85cm,y=.8cm,font=\footnotesize,>=latex]
  \begin{scope}
    \node at (2.5,2.95) {Merging two rectangles};
    \draw[dashed,gray] (0,0) rectangle (5.0,1.95);
    \filldraw[fill=blue!8,draw=blue!55!black] (0,0) rectangle (2.05,1.10);
    \filldraw[fill=blue!8,draw=blue!55!black] (2.95,.65) rectangle (5.0,1.95);
    \foreach \x/\y in {0/.4,.7/1.10,1.45/0,2.05/.70,2.95/1.33,3.65/.65,4.38/1.95,5.0/1.32}
      \fill[blue!65!black] (\x,\y) circle (1.7pt);
    \draw[<->] (2.05,.41) -- (2.95,.41);
    \node[above,inner sep=1pt] at (2.50,.43) {$\le2t$};
    \draw[gray] (0,-.45) -- (5,-.45);
    \foreach \x in {0,.7,1.45,2.05,2.95,3.65,4.38,5}
      \fill[blue!65!black] (\x,-.45) circle (1.5pt);
    \node at (2.5,-.99) {Coordinate projection};
  \end{scope}
  \begin{scope}[xshift=5.95cm]
    \node at (2.45,2.95) {Separated marks on a long side};
    \draw[gray] (-.18,0) rectangle (5.13,1.95);
    \foreach \x/\y in {.55/.53,2.42/1.25,4.45/.64}{
       \draw[blue!62!black,densely dashed,fill=blue!7]
         (\x-.43,\y-.43) rectangle (\x+.43,\y+.43);
       \fill[blue!65!black] (\x,\y) circle (2.1pt);
    }
    \foreach \x/\y in {0/.6,.29/.17,.79/1.09,1.08/.43,1.95/1.53,2.30/.52,2.74/1.86,3.03/.99,3.97/1.07,4.69/.28,4.95/1.34}
        \fill[gray] (\x,\y) circle (1.25pt);
    \draw[<->] (-.18,-.42)--(5.13,-.42);
    \node[below] at (2.47,-.46) {side length at least $L$};
  \end{scope}

\end{tikzpicture}
\caption{Merging preserves the bound $2t$ on projected gaps. A long
side therefore forces many marks at pairwise distance greater than $2R$.
Their radius-$R$ input boxes (dashed) are disjoint.}
\label{fig:slabs}
\end{figure}

\begin{lemma}[One scale]\label{lem:slabs}
Fix a positive integer $t$. Let $Q$ be $R$-local with
$\sup_x\PP(x\in Q)\le e^{-K}$, where $K\ge d^2$, and set
$L=\lfloor e^{K/(16d)}\rfloor$. Suppose $R\le\sqrt L$ and $n>8L$.
For all large $d$ there are rectangles $\mathcal R$ of maximum side
length less than $L$, pairwise at $\ell^\infty$ distance greater than
$2t$, and a $(3L+R)$-local set $Q'$ such that, writing
$U=\bigcup_{D\in\mathcal R}D$,
\[
 Q\subset U\cup Q',\qquad
 \sup_x\PP(x\in Q')\le e^{-\sqrt L},\qquad
 \sup_x\PP(B^\infty(x,t)\cap U\ne\varnothing)\le e^{-K/2}.
\]
\end{lemma}
\begin{proof}
Start with one singleton for each marked site in $Q$. Whenever two
active rectangles are at $\ell^\infty$ distance at most $2t$, replace them by their
coordinatewise hull, breaking ties in a fixed order. Use lifts in which
a closest pair of sites differs by at most $2t$ in each coordinate.
The new sides have length at most $2L+2t<n/2$ for large $d$. Set the
rectangle aside as soon as a side has length at least $L$. A rectangle
set aside takes no further part in the merging. Active rectangles may
later come close to it, which is harmless because its sites will be
covered by $Q'$. The active rectangles left at the end form
$\mathcal R$. Each contains an original mark. Thus $B^\infty(x,t)$ can meet their union only if there is a mark
within distance $L+t$ of $x$, an event of probability at most
$(2L+2t+1)^d e^{-K}\le e^{-K/2}$.

In a lift of a merged rectangle, the projections of its original marks
contain both endpoints of every side, with consecutive projected marks
at most $2t$ apart. This holds for singletons and is preserved when two
projected intervals overlap or have distance at most $2t$
(Figure~\ref{fig:slabs}). Put $b=\lceil L/(4R+2t+1)\rceil$. Declare
$x\in Q'$ if $B^\infty(x,3L)$ contains $b$ marks at pairwise toroidal
$\ell^\infty$ distances greater than $2R$.

In a rectangle set aside, take a maximal such separated family of its
original marks. Every original mark is within $2R$ of a chosen one.
On a longest side the intervals of radius $2R+t$ about the chosen
projections therefore cover the side. These distances can be read in
the lift since each side is shorter than $n/2$. Hence there are at
least $b$ chosen marks, all within $3L$ of every point of the
rectangle. Thus $Q'$ covers the discarded rectangles and has dependence
radius $3L+R$. The $R$-windows of any specified separated family are
disjoint. A union bound gives
\[
 \PP(x\in Q')\le\big[(6L+1)^d e^{-K}\big]^b
 \le e^{-Kb/2}\le e^{-\sqrt L},
\]
using $K\ge d^2$ and $R\le\sqrt L$.
\end{proof}

\subsection{Stopping the process}
To keep track of the side lengths, put $\T_0(x)=x$ and
$\T_{j+1}(x)=2^{\T_j(x)}$. Thus $\T_j(1)=\T_j$.

\begin{proposition}[A local test for nonpercolation]\label{prop:cover}
Fix an integer $h\ge1$ and a constant $A>0$. Use the sets in \eqref{eq:marks}, set
$t=h+1$, and let $H\in\{0,\ldots,d\}$ be a budget for the number of
rectangle scales.
Suppose $3\le n\le\T_H(d^A)$,
$K/d^2\to\infty$, and $K\le2^d$. If every deterministic union $V$ of
at most $H$ proper rectangles satisfies
\begin{equation}\label{eq:late-local}
 \sup_x\PP_p(x\in E_h(V))\le e^{-2K},
\end{equation}
and $\sup_x\PP_p(x\in B(\varnothing))\le3/4$, then for large $d$,
\begin{align}
 \PP_p(A_\infty=\ZZ_n^d)&\le e^{-cK},\label{eq:global-fail}\\
 \PP_p(x\in A_\infty\setminus A_0)
 &\le\PP_p(x\in B(\varnothing)\setminus A_0)+(d+1)e^{-K/2}.
 \label{eq:density-comparison}
\end{align}
\end{proposition}
\begin{proof}
Mark $x$ if $x\in E_h(V)$ for some such $V$. The initial bits and the
trace of $V$ in $B^\infty(x,t)$ determine this event, so the marks
form a $t$-local set $Q_0$. Each coordinate interval has at most
$2^{2t+1}$ traces on its coordinate window, including intervals which
wrap around the torus. Fix one deterministic representative for each
exact trace of a union. There are at most $\exp\{C_t d(H+1)\}$
representatives. Equation~\eqref{eq:late-local} and a union bound give
$\PP_p(x\in Q_0)\le e^{-K}$.

Starting at $K_0=K$, $R_0=t$, define
\[
 L_i=\lfloor e^{K_i/(16d)}\rfloor,\qquad
 K_{i+1}=\sqrt{L_i},\qquad R_{i+1}=3L_i+R_i .
\]
Apply Lemma~\ref{lem:slabs} until the first $J$ with $n\le8L_J$.
There are at most $H$ steps. Indeed, $L_0^{1/4}\ge d^A$ for large $d$,
and
$\log(L_{i+1}^{1/4})\ge\sqrt{L_i}/(64d)-1\ge L_i^{1/4}$.
Thus $L_i^{1/4}\ge\T_i(d^A)$ and $R_i\le\sqrt{L_i}$ by induction.
Since $n\le\T_H(d^A)\le L_H^{1/4}$, the stopping index satisfies $J\le H$.
At the stopping scale,
$\PP_p(Q_J\ne\varnothing)\le (8L_J)^d e^{-K_J}\le e^{-K_J/2}$.

Let $U$ be the union of all retained rectangles. Each has sides shorter
than $n/8$ and is proper. On $B^\infty(x,t)$ at most one rectangle per
scale meets the window, so $U$ agrees there with one of the exact
traces just tested. Hence $x\in E_h(U)$ implies $x\in Q_0$, even
though $U$ was chosen from the seeds. When $Q_J=\varnothing$, the
successive covers give $E_h(U)\subset Q_0\subset U\subset B(U)$.
By~\eqref{eq:marks}, $A_\infty\subset B(U)$.

A fixed window meets $U$ with probability at most $d e^{-K/2}$
by Lemma~\ref{lem:slabs}. If it misses $U$, membership at its centre
in $B(U)$ agrees with membership in $B(\varnothing)$. This proves
\eqref{eq:density-comparison}. For \eqref{eq:global-fail}, choose
$\lceil K\rceil$ vertices at pairwise graph distance greater than $2t$.
A greedy packing has size at least $3^d/(2d+1)^{2t}>2^d+1$.
The empty-background events use disjoint graph balls, so the probability
that all the vertices belong to $B(\varnothing)$ is at most
$(3/4)^{\lceil K\rceil}$. The probability that a selected window meets
$U$, or that $Q_J$ is nonempty, is at most
$(d\lceil K\rceil+1)e^{-K/2}$. Adding these bounds proves the claim.
\end{proof}

\begin{proof}[Proof of Theorem~\ref{thm:density}]
Take $H=H(n)$. The hypothesis gives $D_H\ge\varepsilon d$ for some
$\varepsilon>0$ and all large $d$. Apply Lemma~\ref{lem:three} and
Proposition~\ref{prop:cover} with $h=2$, $A=1$ and $K=c_0d^3$ for a
sufficiently small $c_0>0$. Since $p<1/2$, the empty-background
probability is at most $p+e^{-cd}<3/4$. In particular
\eqref{eq:global-fail} gives $\PP_p(\ZZ_n^d\text{ fills})\le
e^{-c''d^3}$. At each site \eqref{eq:density-comparison} bounds the
probability of a new infection by $e^{-c'''d}$. Averaging and applying
Markov's inequality gives the convergence in probability.

For the box, extend its seeds to a torus of side $\max\{n,3\}$. This
changes the height by at most one. Deleting vertices and edges gives the
same upper bound on the probability of a new infection at each box site.
When $n\ge3$, filling the box also forces the torus to fill. For bounded
$n$, the $2^{d-1}$ disjoint corner edge-lines give a filling probability
at most $\exp\{-2^{d-1}(1-p)^n\}$ by \eqref{eq:empty-lines}.
\end{proof}

\subsection{The first round accounts for the increase}
\begin{proposition}[Late infections]\label{prop:first-round}
Fix $0<p<1/2$ and suppose $\limsup_{d\to\infty} H(n)/d<\rho_p$.
Put $q=1-p$. On the torus let $\tau=|N(x)\cap N(y)|$ for any
edge $xy$. Then
\begin{align}
 \PP_p(x\in A_\infty\setminus A_1)
    &=(1+o(1))\,\frac{q^2}{2\pi p(q-p)}(4pq)^{\,2d-\tau},
       \label{eq:two-site}\\
 \frac{\EE_p|A_\infty\setminus A_0|}{n^d}
    &=q\,\PP\{\Bin(2d,p)\ge d\}\bigl(1+O(e^{-cd})\bigr).\notag
\end{align}
On every box $[n]^d$ with $n\ge2$ in the same height range, including
at its boundary, one has
\begin{equation}\label{eq:box-density}
 \PP_p(x\in A_\infty\setminus A_1)
    \le e^{-cd}\,\PP_p(x\in A_1\setminus A_0).
\end{equation}
Thus the expected first-round count also gives the total expected
increase on the box up to a factor $1+O(e^{-cd})$. Here $c>0$ may
depend on $p$ and the gap in the height condition.
\end{proposition}
\begin{proof}
Write $I_p=-\log(4pq)>0$. For a box of side $n$, extend its seeds
to the torus of side $\max\{n,3\}$ as in the proof of
Theorem~\ref{thm:density}. In the torus construction put
$E=\{v:d_{A_0}(v)\ge d-G\}$ and $a=h-1$.
Here $h=2$ and $G=2g$. For any fixed radius $R$, except on an event
of probability at most $e^{-K/3}$, $Q_J$ is empty and all the windows
about vertices in $B(x,R)$ miss $U$. Every new infection there, on
either graph, then lies within torus distance $a$ of $E$. Indeed it
belongs to $B(\varnothing)\setminus A_0$ in the dominating torus
construction, and tracing the sets $E_i(\varnothing)$ gives the claim.

Fix an integer $M\ge2$ and take $R=(2a+1)M$. Write
$\lambda=\log(q/p)$. The probability that $B(x,R+a)$ contains more
than $M$ sites of $E$ is at most
\begin{equation}\label{eq:many-centres}
 d^{O_{M,a}(1)}\exp\{(M+1)[-I_pd+(G+2M)\lambda]\}.
\end{equation}
Indeed any two sites have at most two common neighbours. For a
specified tuple delete all shared bits, at most $2M$ per site.
The remaining counts are independent. The exponential Markov bound
at $\lambda$ gives $\PP(\Bin(2d,p)\ge d-u)\le e^{-I_pd+u\lambda}$.

Off these two exceptional events the connected component of new
infections containing $x$, if present, lies inside $B(x,R)$. A path
reaching distance $R$ would up to that point be covered by at most $M$
radius-$a$ balls with centres in $E$. Along a connected path successive
centres are at distance at most $2a+1$, so its diameter would be at
most $(2a+1)(M-1)+2a=R-1$. Moreover all additions in this component
belong to $E$ when $G\ge M$. Its first addition outside $E$ would
have fewer than $d-G$ initially infected neighbours and could gain
at most $M$ more from previous additions. In particular the component
contains at most $M$ sites. Box paths and box updates use subsets of
the torus edges, so the argument applies to both graphs.

Take $M=2$ on the torus. Choose the fixed $\varepsilon$ in
Lemma~\ref{lem:three} small enough that \eqref{eq:many-centres} is
$O(e^{-(2I_p+c_1)d})$. If $x$ is infected after round one, its component
must then consist of $x$ and a neighbour $y$ infected in round one.
Thus $x,y$ start healthy, $d_{A_0}(x)=d-1$ and $d_{A_0}(y)\ge d$.
Conversely each such event infects $x$ in round two. Two different
choices of $y$ require three sites of $E$. After summing over the
$O(d^2)$ pairs, the error in adding the single-neighbour probabilities
is still $O(e^{-(2I_p+c_2)d})$.

Adjacent torus sites have $\tau$ common neighbours. Condition on $x,y$
being healthy and on the number $Z\sim\Bin(\tau,p)$ of infected common
neighbours. The remaining counts are independent $\Bin(2d-1-\tau,p)$.
Put $P_d=\binom{2d}{d}(pq)^d$. For fixed integers $b,j$, the mass at
$d-j$ of $\Bin(2d-b,p)$ is asymptotic to
$P_d\,2^{-b}p^{-j}q^{j-b}$, by taking the ratio of the binomial
coefficients. Its upper tail there is asymptotic to that mass times
$q/(q-p)$, since the first finitely many successive mass ratios tend
to $p/q<1$ and the ratios decrease. Use $b=1+\tau$, with $j=1+Z$ for
the mass at $x$ and $j=Z$
for the tail at $y$. Multiplying by $2dq^2$, using
$\EE[(q/p)^{2Z}]=(q/p)^\tau$ and
$P_d\sim(4pq)^d/\sqrt{\pi d}$ gives~\eqref{eq:two-site}.
The first-round probability is $q\PP(\Bin(2d,p)\ge d)
=\exp\{-I_pd+O(\log d)\}$, which gives the density formula.

For the box choose $M$ with $(M+1)I_p/2>|\log p|+1$, then choose
$\varepsilon$ so small that $(G+2M)\lambda\le I_pd/2$ for all large
$d$. The total exceptional probability is at most $p^d e^{-c_3d}$.
A late infection at $x$ requires an earlier neighbouring addition $y$,
and both have at least $d-M$ initially infected \emph{box} neighbours.
Their neighbour sets are disjoint by bipartiteness. For $N=\deg(x)\ge d$,
ratios of consecutive binomial masses give
\[
 \PP(\Bin(N,p)\ge d-M)
 \le O_{p,M}(d^M)\,\PP(\Bin(N,p)\ge d).
\]
At $y$ the probability is at most $e^{-I_pd+M\lambda}$. Sum over its
at most $2d$ choices. Since $\PP(\Bin(N,p)\ge d)\ge p^d$ and $q$
is fixed, this proves \eqref{eq:box-density}.
\end{proof}

\section{Growing a box one face at a time}\label{sec:box}
Each vertex on a new face has an infected neighbour behind it. The
initial state of the neighbour ahead of it supplies one more possible
infection. We keep track of both contributions as the box grows.

\subsection{Exterior seeds and the boundary of a box}
For a box $F\subset\ZZ^k$, give each site in $F$ an initial seed bit
and $\ell$ additional \emph{private bits}. Give each vertex outside $F$
adjacent to $F$ a seed bit as well. All these bits are independent
Bernoulli($p$). Only sites in $F$ update. Writing $\beta_x$ for the sum of
the private bits at $x$, the condition for infection is
\begin{equation}\label{eq:face-rule}
 \#\{\text{infected neighbours in }F\}
 +\#\{\text{exterior neighbours that are seeds}\}+\beta_x\ge k .
\end{equation}
Denote this process by $\mathcal M(k,\ell)$. Every site, including a
corner, has $2k$ neighbour bits. If all sides have length at least two,
at least $k$ of its neighbours belong to $F$ and can become infected.
On a torus there are no exterior sites. For each integer $m\ge1$, let
$g_{k,\ell}(m)$ be the maximum failure probability over boxes with
\emph{every} side length in $[m,2m]$. Set $g_{0,\ell}(m)=0$, since
the threshold on a zero-dimensional cell is zero.

\begin{figure}[htbp]
\centering
\begin{tikzpicture}[x=5.8mm,y=5.8mm,font=\footnotesize,>=latex]
  \begin{scope}
    \fill[blue!13] (.5,.5) rectangle (5.5,4.5);
    \foreach \x in {0,...,6}{\draw[black!25] (\x,0) -- (\x,5);}
    \foreach \y in {0,...,5}{\draw[black!25] (0,\y) -- (6,\y);}
    \foreach \x in {1,...,5}{\foreach \y in {1,...,4}{\fill[blue!58!black] (\x,\y) circle (1.55pt);}}
    \foreach \x in {1,...,5}{\foreach \y in {0,5}{\fill[orange!78!black] (\x,\y) circle (1.85pt);}}
    \foreach \x in {0,6}{\foreach \y in {1,...,4}{\fill[orange!78!black] (\x,\y) circle (1.85pt);}}
    \foreach \x in {0,6}{\foreach \y in {0,5}{\draw[fill=white,thick] (\x,\y) circle (2.1pt);}}
    \node[below] at (3,-.35) {Interior, then four sides, then four corners};
  \end{scope}
  \begin{scope}[xshift=7.0cm]
    \fill[blue!13] (-.45,-.45) rectangle (4.45,2.5);
    \foreach \x in {0,...,4}{\draw[black!25] (\x,-.3) -- (\x,3);}
    \foreach \y in {0,...,3}{\draw[black!25] (0,\y) -- (4,\y);}
    \draw[orange!76!black,line width=1.5pt] (0,3) -- (4,3);
    \draw[black!55] (-1,3) -- (0,3);
    \filldraw[fill=white] (-1,3) circle (2.7pt);
    \fill[orange!78!black] (0,3) circle (2.8pt);
    \fill[blue!58!black] (0,2) circle (2.8pt);
    \fill[orange!78!black] (1,3) circle (2.8pt);
    \node[above=3pt] at (-1,3) {seed};
    \node[above=3pt] at (.12,3) {$x$};
    \draw[->,blue!50!black] (1.2,1.2) -- (.15,1.9);
    \node[below right,inner sep=1pt,fill=blue!13] at (.55,1.18) {inward support};
    \node[below] at (1.75,-.9) {An endpoint of a side cell};
  \end{scope}
\end{tikzpicture}
\caption{The interior fills before the sides and corners. At the endpoint
$x$ of the top side cell, the neighbour below is infected and the corner
to the left supplies its initial state.}
\label{fig:cells}
\end{figure}
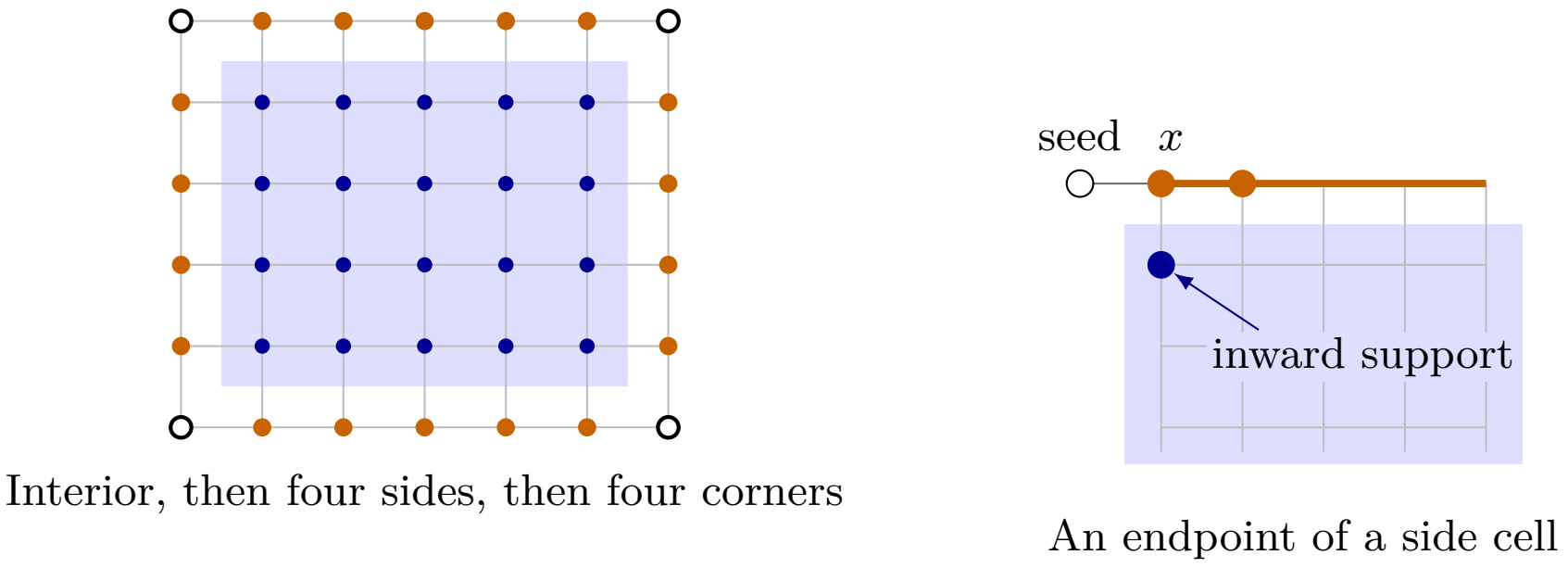

\begin{lemma}[The boundary faces]\label{lem:cells}
For the ordinary $d$-neighbour process on $[n]^d$, $n\ge4$,
\[
 \PP_p([n]^d\text{ fails to fill})
 \le \sum_{j=1}^d \binom{d}{j}2^{d-j}\,g_{j,0}(n-2).
\]
\end{lemma}
\begin{proof}
Split each coordinate interval into $\{1\},[2,n-1],\{n\}$. Fill the
product cells in increasing order of codimension, starting with the
interior. A cell with $d-j$ fixed coordinates then has one infected
inward neighbour in each fixed direction. Both neighbours in every
free coordinate belong to the original box, even at an endpoint of the
cell (Figure~\ref{fig:cells}). Test this cell using its own seeds and
the seeds just beyond its ends. Prescribing the $d-j$ inward neighbours
as infected leaves the rule~\eqref{eq:face-rule} with $k=j$ and no
private bits, that is, $\mathcal M(j,0)$. Each test uses the original
bits, even where different cells share a seed. If all cell tests succeed,
the cells fill in the stated order and the corners fill from their $d$
inward neighbours. There are $\binom dj2^{d-j}$ cells of dimension $j$.
\end{proof}

\subsection{A path through healthy sites}
A nonempty component of the final healthy set must meet opposite faces
in some coordinate. Otherwise choose a missed face in each coordinate
and orient that coordinate towards it. At a vertex of the component
maximizing the sum of the oriented coordinates, all $k$ neighbours in
the chosen directions belong to the box and lie outside the component.
They are infected and force the vertex to become infected.

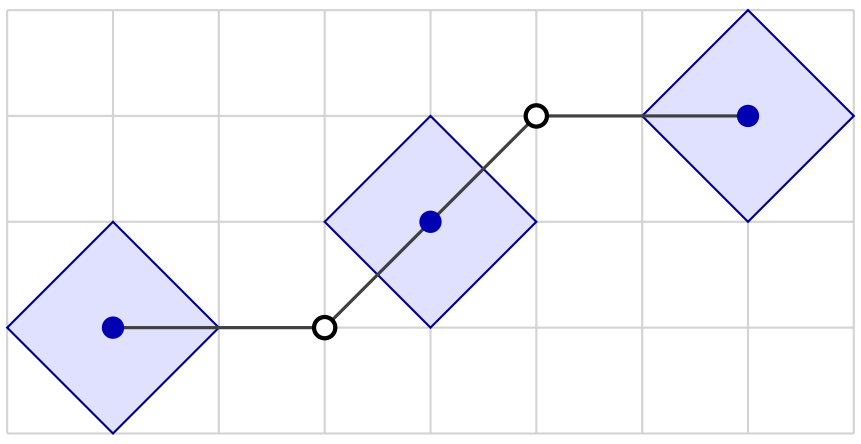
\begin{figure}[htbp]
\centering
\begin{tikzpicture}[x=7mm,y=7mm,font=\footnotesize]
  \draw[step=1,black!17] (-1,-1) grid (7,3);
  \foreach \x/\y in {0/0,3/1,6/2}{
    \filldraw[fill=blue!12,draw=blue!54!black,thin]
      (\x-1,\y) -- (\x,\y+1) -- (\x+1,\y) -- (\x,\y-1) -- cycle;
  }
  \draw[black!75,semithick] (0,0) -- (2,0) -- (3,1) -- (4,2) -- (6,2);
  \foreach \x/\y in {0/0,3/1,6/2}{\fill[blue!70!black] (\x,\y) circle(2.1pt);}
  \foreach \x/\y in {2/0,4/2}{\draw[fill=white,thick] (\x,\y) circle(2pt);}
\end{tikzpicture}
\caption{A shortest path with steps of lattice distance at most two.
The shaded radius-one neighbourhoods at alternate vertices are disjoint.}
\label{fig:path}
\end{figure}

\begin{lemma}[A short path test]\label{lem:path}
Let $S$ be a random subset of the vertices of a finite graph $F$ of maximum degree $D$, with
$\PP(x\in S)\le e^{-Q}$ at each site. For an integer $t\ge1$, assume that the
membership events are independent on every set of sites at pairwise
graph distance greater than $2t$. If $Q\ge8t\log(D+1)$, then the
probability of a nearest-neighbour component of $S$ containing two
sites at distance at least $a$ is at most
\[
 |F|\exp\{-Qa/(8t)\}.
\]
\end{lemma}
\begin{proof}
Join sites of $S$ whenever their distance is at most $2t$. A shortest
path between the two sites in this graph has at least $a/(2t)$ steps.
Nonconsecutive vertices of a shortest path are more than $2t$ apart,
so membership at the even-indexed vertices is independent
(Figure~\ref{fig:path}). Count initial
segments of $u=\lceil a/(2t)\rceil$ steps having this separation.
There are at most $|F|(D+1)^{2tu}$ of them, and each has probability
at most $e^{-Qu/2}$. This proves the estimate.
\end{proof}

On a $k$-dimensional torus at threshold $k$, each vertex of a final
healthy component $C$ has at least $k+1$ neighbours in $C$, even with
private bits. Fix $x\in C$ and let $v\in C$ maximize
$i=\operatorname{dist}(x,v)$ in the ambient graph. Its neighbours
in $C$ all belong to $B(x,i)$, so \eqref{eq:backward} gives $i\ge k/2$.
Some coordinate projection of $C$ must also be the
whole cycle. Otherwise cut each cycle at a missed level. At a maximizer
of the sum of lifted coordinates, all $k$ positive neighbours lie
outside $C$ and hence are infected. Thus $C$ has diameter at least
$\lfloor n/2\rfloor$ as well. The time-$t$ survivor indicators use
disjoint input balls at pairwise distances greater than $2t$.
If their one-site probabilities are at most $e^{-Q}$ and
$Q\ge8t\log(2k+1)$, Lemma~\ref{lem:path} at threshold $k$ gives
\begin{equation}\label{eq:torus-path}
 \PP_p(\ZZ_n^k\text{ fails})\le
 n^k\exp\!\left\{-\frac{Q}{8t}\max\{k/2,\lfloor n/2\rfloor\}\right\}.
\end{equation}
In particular the bound is at most $e^{-k}$ uniformly in $n\ge3$ when
$Q/t\gg\log k$, since $k\log n/\max\{k,n\}\le\log k$ for $k,n\ge3$.

Fix $k+\ell\le D$ and put $\delta=e^{-100D}$. For all large $D$,
\begin{equation}\label{eq:block-tail}
 g_{k,\ell}(m)\le\delta
 \quad\Longrightarrow\quad
 g_{k,\ell}(n)\le e^{-5Dn/m}\qquad(n\ge m).
\end{equation}
For $k\ge1$, divide a side of length $a\in[n,2n]$ into $\lfloor a/m\rfloor$
nearly equal intervals. The resulting blocks have side lengths in
$[m,2m]$. Test each block with its own private bits and the initial
states of its exterior neighbours. A successful test fills the block in
the larger process. Blocks at pairwise graph distance greater than two
in the index lattice have disjoint inputs. A surviving component gives
a connected set of failed blocks that meets opposite faces. Put
$\xi=n/m$. If $\xi\ge8$, it spans a block-index distance at least
$\lfloor\xi\rfloor-1\ge3\xi/4$. Lemma~\ref{lem:path} with $t=1$
gives a failure bound $(2\xi)^k e^{-100D(3\xi/4)/8}\le e^{-5D\xi}$.
For $1\le\xi<8$, a union bound over at most $16^k$ blocks suffices.

\subsection{One exponential for a new dimension}
\begin{proposition}[Growth from smaller faces]\label{prop:growth}
Fix $0<p_-<1$, take $p\in[p_-,1)$, and put $\delta=e^{-100D}$.
Among the nonnegative integer pairs $k+\ell\le D$, designate a set
of base pairs containing every pair with $k\le1$. Suppose that, for
an integer $S\ge1000D$, each base pair satisfies $g_{k,\ell}(S)\le\delta$.
Let $h=h(k,\ell)$ be the number of steps
$(k,\ell)\mapsto(k-1,\ell+1)$ to the first base pair, with $h=0$
at a base pair. Put
\[
 R_0=S,\qquad R_{i+1}=\lceil\exp(R_i^{\,4D})\rceil .
\]
Then, for all sufficiently large $D$ depending only on $p_-$,
\[
 g_{k,\ell}(R_h)\le\delta,\qquad
 g_{k,\ell}(n)\le e^{-5Dn/R_h}\quad(n\ge R_h).
\]
If $S$ is bounded by a fixed polynomial in $D$, then, for all large
$D$, $R_i\le \T_{H(D)+2+i}^{1/(8D)}$ for every $i\ge0$.
\end{proposition}

\begin{proof}
It suffices by~\eqref{eq:block-tail} to prove the first bound. Induct
on $h$. The case $h=0$ is the hypothesis. Otherwise let
$R=R_{h-1}$, $M=R^2$ and $N=R_h$. Fix a cube $Q$ of side $M$ in a
target box of side lengths $L_i\in[N,2N]$.

Starting from $Q$, visit the coordinates cyclically. At each visit add
one site on the right if possible and otherwise on the left. Skip
intervals already full. After $\tau-M$ cycles the lengths are
$\min(\tau,L_i)$. Test each new face with its inward neighbours prescribed
infected. The outward neighbour of each site supplies one extra private
bit, even at the target boundary, where it is an exterior seed.
Neighbours just beyond the edges of the face supply its in-plane
exterior seeds. The outward bits are distinct and lie outside the
face plane, so they are independent of all the seed and private bits
used in that plane. Thus each test has law $\mathcal M(k-1,\ell+1)$.
Define all tests on the original inputs. If $Q$ is full and the tests
succeed, the target fills (Figure~\ref{fig:face}).

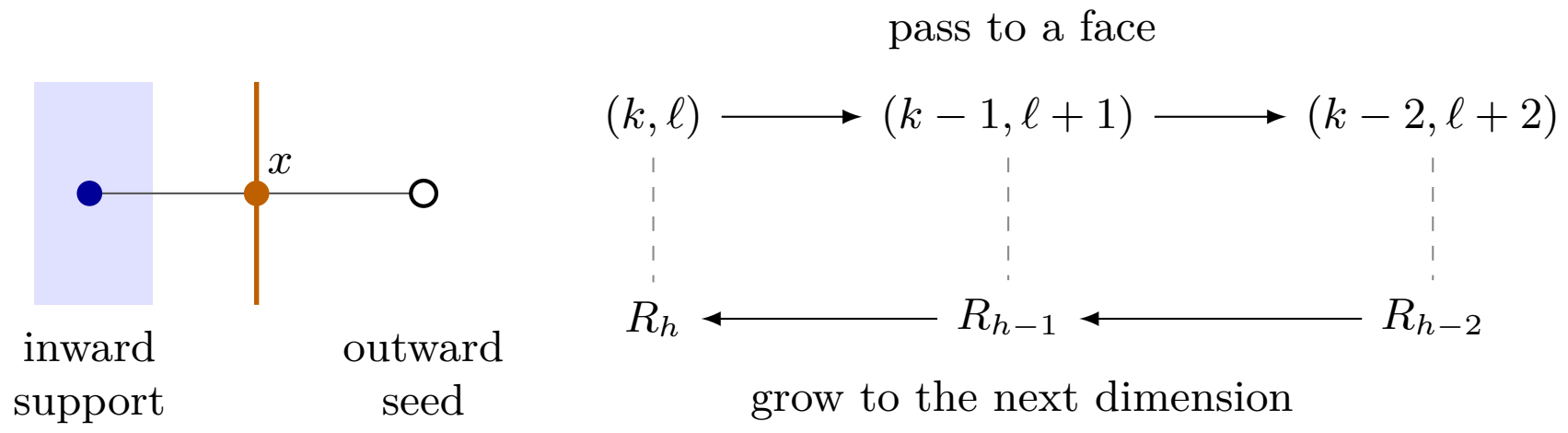
\begin{figure}[htbp]
\centering
\begin{tikzpicture}[x=1cm,y=1cm,font=\footnotesize,>=latex]
  \begin{scope}
    \fill[blue!12] (-.4,-.35) rectangle (.45,1.25);
    \draw[orange!74!black,line width=1pt] (1.2,-.35)--(1.2,1.25);
    \draw[black!65] (0,.45)--(2.4,.45);
    \fill[blue!60!black] (0,.45) circle (2.6pt);
    \fill[orange!75!black] (1.2,.45) circle (2.6pt);
    \draw[fill=white,thick] (2.4,.45) circle (2.6pt);
    \node[above right,inner sep=2pt] at (1.2,.52) {$x$};
    \node[below,align=center] at (0,-.40) {inward\\support};
    \node[below,align=center] at (2.4,-.40) {outward\\seed};
  \end{scope}
  \begin{scope}[xshift=4.05cm]
    \node (A) at (0,1) {$(k,\ell)$};
    \node (B) at (2.55,1) {$(k-1,\ell+1)$};
    \node (C) at (5.6,1) {$(k-2,\ell+2)$};
    \node (a) at (0,-.45) {$R_h$};
    \node (b) at (2.55,-.45) {$R_{h-1}$};
    \node (c) at (5.6,-.45) {$R_{h-2}$};
    \draw[->] (A)--(B); \draw[->] (B)--(C);
    \draw[->] (c)--(b); \draw[->] (b)--(a);
    \foreach \typenode/\scalenode in {A/a,B/b,C/c}
      {\draw[black!45,dashed] (\typenode.south)--(\scalenode.north);}
    \node[above] at (2.65,1.35) {pass to a face};
    \node[below] at (2.65,-.76) {grow to the next dimension};
  \end{scope}
\end{tikzpicture}
\caption{At $x$, the inward support lowers the threshold and the outward
seed adds a private bit. Each step back from a smaller face costs one
exponential in side length, as shown by the lower arrows.}
\label{fig:face}
\end{figure}

Let $G_Q$ be the event that all tests in this growth from $Q$ succeed.
There are at most $k$ faces per cycle. For $M\le\tau<N$, every face added during the cycle from $\tau$ to
$\tau+1$ has side lengths in $[\tau,\tau+1]$. All faces added in later cycles have side lengths in $[N,2N]$. The induction hypothesis
and a geometric sum give
\[
 \PP_p(G_Q^c)\le k(R+1)e^{-M/R}+2Nk e^{-N/R}\le e^{-R/2}
\]
for large $D$, since $M=R^2$ and $R\ge1000D$. Each face and its
in-plane exterior sites lie in a plane strictly beyond $Q$, and the
outward sites lie one step farther away. Thus $G_Q$ is independent of
the initial fullness of $Q$, whose probability is $v=p^{M^k}$.

The target contains at least $(N/(2M))^k$ disjoint candidate cubes.
This exceeds $200D/v+1$ for large $D$, since $M^k\le R^{2D}$,
$\log N\ge R^{4D}$ and $p\ge p_-$. Retain $b=\lceil200D/v\rceil$ of
them. If the target fails, either none of these cubes is full or a full
one fails its own growth test. Independence of the initial fullness
events and a union bound give
\[
 \PP_p(\text{failure})\le(1-v)^b+bv e^{-R/2}
 \le e^{-200D}+(200D+1)e^{-R/2}\le\delta .
\]

For the height estimate let $Y_0=\T_{H(D)+2}\ge2^{2^D}$. Then
$S^{8D}\le Y_0$ for large $D$. Whenever $R_i\le Y^{1/(8D)}$ with
$Y\ge Y_0$,
\[
 R_{i+1}\le e^{2\sqrt Y}\le (2^Y)^{1/(8D)}.
\]
Induction gives the assertion.
\end{proof}

\subsection{Starting at a linear side length}
Write
\[
 s=k+\ell,\qquad q=1-p,\qquad
 \Delta=k-p(2k+\ell)=(1-2p)s-q\ell .
\]
Thus $\Delta$ is the threshold minus the mean initial support. Passing
to a face preserves $s$ and reduces $\Delta$ by $q$.

\begin{lemma}[Failure on a fixed-dimensional face]\label{lem:fixed-face}
Fix $k\ge1$, $\ell\ge0$ and $p\in(0,1)$. With $\gamma=-\log(1-p)$,
\[
 \lim_{m\to\infty}-\frac1m\log g_{k,\ell}(m)=(k+\ell)\gamma .
\]
\end{lemma}
\begin{proof}
In $[m]^k$ keep a specified corner edge-line healthy. It suffices
that every site bit on the line, its private bits, its $k-1$
transverse exterior bits per site, and its two axial exterior bits are
zero. The other $k-1$ neighbours at each site cannot infect it. The
probability is $(1-p)^{(k+\ell)m+2}$, exactly the failure probability
on the line when $k=1$. For an upper bound on failure it suffices to
treat $[m]^k$. Its $2^k$ translates flush with the corners cover
any box with sides in $[m,2m]$. If their individual tests on the
original bits pass, the box fills.

Induct on $k$, keeping $s=k+\ell$ fixed. Write $b=\lfloor m^{1/3}\rfloor$.
In a slab $[m]^{k-1}\times[b]$ place $\lfloor m/b\rfloor$ disjoint
$b$-cubes along one long coordinate. For fixed $k,p$ each spans from
its own seeds at threshold $k$ with probability at least $1/2$ for
large $b$~\cite[Theorem~1]{BBDCM}. An exponential bound for this
fixed-dimensional estimate is given in~\cite[\S2.11, Lemma~2.45]{Book}.
Test every coordinate $(k-1)$-cube of side $b$ in the slab in both
growth directions. Prescribe infected inward neighbours and use outward
seed bits as one additional private bit per site, as in
Proposition~\ref{prop:growth}. If a seed cube spans and every test
passes, grow one layer at a time in the long coordinates. Each layer
is covered by tested cubes with infected neighbours behind them.
Each test has law $\mathcal M(k-1,\ell+1)$, so the probability of
slab failure is at most
\[
 2^{-\lfloor m/b\rfloor}+2km^k g_{k-1,\ell+1}(b)
 \le \exp\{-(s\gamma-o(1))b\}.
\]
Here $\log m=o(b)$ and $b=o(m/b)$. Place $\lfloor m/(b+2)\rfloor$
such slabs in the target, with two unused layers between successive
slabs. Their tests use disjoint inputs, including the one-layer
exterior boundaries. Thus the probability all fail is at most
$\exp\{-(s\gamma-o(1))m\}$. A successful slab grows through all
remaining layers unless a full-width $(k-1)$-face test fails in one
of the two orientations. These tests cost at most
$2m g_{k-1,\ell+1}(m)$, which proves the induction.
\end{proof}

\begin{lemma}[Starting from a local estimate]\label{lem:base}
Fix $0<p_-<1$, take $p\in[p_-,1)$, put $\alpha=p_-^2/100$, and fix an
integer $s_0\ge2$.
Let $W(s)/s$ be nonincreasing for $s>0$. Designate as base pairs those
with $k\le1$, $s\le s_0$, or $\Delta\le W(s)$. For every base pair
with $s>s_0$ and $k\ge\max\{2,\alpha s\}$, suppose there are an integer
$t\ge1$ and a number $Q$ such that, in every box of side lengths at
least two with exterior seeds, all sites satisfy
\begin{equation}\label{eq:start-local}
 \PP_p(x\text{ is healthy at time }t)\le e^{-Q},\qquad
 Q\ge16t\bigl(1+\log(2k+1)\bigr).
\end{equation}
Then $g_{k,\ell}(\lceil CD\rceil)\le e^{-100D}$ at every base pair,
where $C=C(p_-,s_0)$, so Proposition~\ref{prop:growth} applies with
$S=\lceil CD\rceil$. For the ordinary box $[n]^d$ take $D=d$ and put
\begin{equation}\label{eq:potential}
 J=\max\left\{0,\left\lceil\frac{(1-2p)d-W(d)}q\right\rceil\right\}.
\end{equation}
For all sufficiently large $d$, the box fills with probability at least
$1-e^{-d}$ whenever $n\ge R_J+2$.
\end{lemma}
\begin{proof}
Write $m=\lceil CD\rceil$. A line containing any initial seed fills, so
$g_{1,\ell}(m)\le e^{-p_-m}$. For $s\le s_0$, Lemma~\ref{lem:fixed-face} at $p_-$ gives
$g_{k,\ell}(m)\le e^{-c_0m}$ for large $m$, uniformly over these finitely
many pairs and $p\ge p_-$.

If $s>s_0$ and $2\le k<\alpha s$, the private bits alone have mean
$p\ell\ge k+p_-s/2$. Hoeffding's inequality gives~\eqref{eq:start-local}
with $t=1$ and $Q=p_-^2s/2\ge50k$. At the remaining base pairs use the
hypothesis. A component that survives in a box with side lengths in
$[m,2m]$ crosses it, and all its sites were healthy at time $t$.
Those states depend on inputs within distance $t$, including any
exterior seeds. Lemma~\ref{lem:path} bounds the crossing probability by
\[
 (2m)^k \exp\{-Q(m-1)/(8t)\}\le e^{-m}.
\]
Here $k\log(2m)\le m\log(2k)$ for $2\le k\le m$, since
$\log(2x)/x$ decreases on $[2,\infty)$. Taking $C$ large makes all
the base probabilities at most $e^{-100D}$.

For a $j$-dimensional cell in Lemma~\ref{lem:cells}, the pair
$(j,0)$ has deficit $(1-2p)j$. After $i$ steps the deficit is
$(1-2p)j-qi$. The number of steps to a base pair is at most the
positive part of $\lceil((1-2p)j-W(j))/q\rceil$. This is at most $J$,
since $j\max\{0,1-2p-W(j)/j\}$ is nondecreasing in $j$. Thus for
$n-2\ge R_J$, the probability that any cell fails is at most
$3^d e^{-5d}\le e^{-d}$.
\end{proof}

\subsection{The leading threshold}

\begin{proposition}[Filling above the curve]\label{prop:leading-up}
Fix $0<p\le1/2$ and $\varepsilon>0$. The box $[n]^d$ fills with
probability at least $1-e^{-d}$, for all large $d$, whenever
\[
 n\ge\T_{\lceil(1-2p+\varepsilon)d/(1-p)\rceil+H(d)+2}.
\]
At $p=1/2$ the same probability bound holds on boxes with $n\ge Cd$,
for an absolute $C$, and on every torus with $n\ge3$.
\end{proposition}
Proposition~\ref{prop:fine-up} sharpens this bound. We keep the simpler
argument, so that Theorem~\ref{thm:phase} does not depend on
Section~\ref{sec:fine}.
\begin{proof}
Take $W(s)=-\varepsilon s$ in Lemma~\ref{lem:base}. At a base pair with
$s>s_0$ and $k\ge\max\{2,\alpha s\}$, the mean initial support is at
least $k+\varepsilon s$ at every site in the process with exterior
seeds. Its chance of staying healthy through the first round is at most
$e^{-\varepsilon^2s}$.
This gives~\eqref{eq:start-local} with $t=1$ on taking $s_0$ large.
The height $J$ in~\eqref{eq:potential} is
$\lceil(1-2p+\varepsilon)d/(1-p)\rceil$. At the stated side length,
Proposition~\ref{prop:growth} gives $R_J+2\le n$.

At $p=1/2$ take $W=0$. Since $\Delta=-\ell/2$, all pairs are base.
In a box with exterior seeds (all sides at least two) or on a torus,
a vertex $x$ has at least $k$ neighbours that update. Count those
initially healthy neighbours with at least $k$ seeds among their $2k-1$
neighbours other than $x$. By binomial symmetry each qualifies with
probability $1/4$. A seed other than $x$ enters at most four such
neighbour counts by~\eqref{eq:backward}, and its own initial state
enters at most one further test. Thus~\eqref{eq:read} with $b=5$ bounds
the probability of fewer than $k/8$ qualifying neighbours by
$e^{-k/160}$. The initial support at $x$ has mean at least $k$, so
Hoeffding bounds the probability of support below $7k/8$ by
$e^{-k^2/(64s)}$. When both counts meet their bounds, at least $k/8$ new
infections add to the initial support and $x$ is infected by round two. For $k\ge\alpha s$, the sum of the two errors
verifies~\eqref{eq:start-local} with $t=2$ for large $s_0$. Hence
Lemma~\ref{lem:base} gives the box assertion.

On the torus the same two-round calculation has $s=k=d$ and gives
$Q\ge cd$. By~\eqref{eq:torus-path} the probability of failure is at
most $e^{-d}$ for every $n\ge3$.
\end{proof}

\subsection{Boxes of linear side length}\label{sec:short}
Consider the $2^{d-1}$ lines obtained by letting the first coordinate
vary and fixing every other coordinate at $1$ or $n$. On each such
line a site has at most $d-1$ neighbours off the line, so an initially
empty line stays healthy. These lines are disjoint. Thus, for $n\ge2$,
\begin{equation}\label{eq:empty-lines}
 \PP_p([n]^d\text{ fills})\le (1-(1-p)^n)^{2^{d-1}}
 \le \exp\{-2^{d-1}(1-p)^n\}.
\end{equation}
In particular $p_c([n]^d,d)\to1$ for $2\le n=o(d)$.

Suppose $n/d\to a\in(0,\infty)$. Equation~\eqref{eq:empty-lines}
gives $\liminf p_c\ge1-2^{-1/a}$. Theorem~\ref{thm:density} gives
$\liminf p_c\ge1/2$ since $H(n)=o(d)$. For the opposite bound fix $p$
strictly above both limits. Put $m=n-2\sim ad$ and $\gamma=-\log(1-p)$.
For fixed $K$, Lemma~\ref{lem:fixed-face} gives
\[
 \sum_{1\le k<K}\binom dk 2^{d-k}g_{k,0}(m)
 \le 2^d O_K(d^K)e^{-(\gamma-o(1))m}=o(1),
\]
since $a\gamma>\log2$. For large $k$ the two-round estimate at $p=1/2$
in Proposition~\ref{prop:leading-up} gives one-site survival at most
$e^{-ck}$, also an upper bound at $p$. By Lemma~\ref{lem:path}, failure
of a particular $k$-cell has probability at most
$m^k e^{-ck(m-1)/16}\le e^{-c'km}$.
Choose $K$ with $c'Ka>2\log3$ and sum this bound over at most $3^d$
cells. Lemma~\ref{lem:cells} proves filling with high probability and
establishes~\eqref{eq:linear-box}.

\begin{proposition}[Empty lines at the threshold]\label{prop:poisson}
Suppose $n=d$, or $n/d\to a$ for fixed $0<a<1$. Write $q=1-p$ and
$N=d2^{d-1}$. If $Nq^n\to u\in(0,\infty)$, the number of empty
corner edge-lines converges in law to $\operatorname{Poisson}(u)$
and the box fills with probability tending to $e^{-u}$. In particular,
\begin{equation}\label{eq:box-endpoint}
 p_c([d]^d,d)=\frac12+\frac{\log d-\log(2\log2)+o(1)}{2d}.
\end{equation}
\end{proposition}
\begin{proof}
At $n=d$ we have $p\to1/2$ from above, since $N2^{-d}=d/2$.
When $n/d\to a<1$ the limit is
$1-2^{-1/a}>1/2$. The large-dimensional cells pass with high probability
by the preceding two-round bound. Setting $a=1$ in the endpoint case,
choose $p_0$ strictly below the limiting density such that
$2a[-\log(1-p_0)]>\log2$. Lemma~\ref{lem:fixed-face} at $p_0$ and the
same cell count show that all cells of dimensions $2\le k<K$ pass
with high probability. A one-dimensional cell fills exactly when its
line, including its two corner endpoints, contains a seed. Hence the
filling event differs from the absence of empty lines by probability $o(1)$.

Expose the corner bits. The number $Y$ of edge-lines with both endpoints
healthy has mean $Nq^2$ and variance $O(dN)$, since each line shares
an endpoint with at most $2(d-1)$ others. Thus $Y/(Nq^2)\to1$ in
probability. Line interiors are disjoint, so conditionally on the corners
the number of empty lines has law $\Bin(Y,q^{n-2})$. Its success
probability tends to zero and its conditional mean tends in probability
to $u$. The binomial Poisson limit gives both assertions. For $n=d$,
bracket the median using $Nq^d=\log2\pm\eta$, then let $\eta\downarrow0$.
This gives $p_c=1-(\log2/N)^{1/d}+o(d^{-1})$, which expands
to~\eqref{eq:box-endpoint}.
\end{proof}

\begin{proof}[Proof of Theorem~\ref{thm:phase} and Corollary~\ref{cor:first-crossing}]
The linear-side limit \eqref{eq:linear-box} was proved above. For the
tower curve pass to a subsequence on which
$\min\{H(n)/d,1\}\to\rho\in[0,1]$. At fixed $0<p<p_\rho$ we have
$H(n)/d\to\rho<\rho_p$. Theorem~\ref{thm:density} gives failure with
high probability on both graphs, so $\liminf p_c\ge p_\rho$. At
$\rho=1$ this follows from $p_c\ge0$.

Proposition~\ref{prop:leading-up} at $p=1/2$ gives $\limsup p_c\le1/2$
on the torus and on boxes with $n\ge C_0d$ for its absolute $C_0$.
For the remaining boxes with $d\le n<C_0d$, every subsequence has a
further subsequence with $n/d\to a\in[1,C_0]$, and \eqref{eq:linear-box}
gives the same bound. Whenever $p_\rho<p<1/2$, we can choose
$\varepsilon>0$ with $(1-2p+\varepsilon)/(1-p)<\rho$. The height in
Proposition~\ref{prop:leading-up} is then at most $H(n)-1$ for large $d$.
Since $n>\T_{H(n)-1}$, both graphs fill with high probability at $p$.
Together these upper bounds give $\limsup p_c\le p_\rho$.
Compactness and continuity of $p_\rho$ give the uniform estimate.

For the first crossing fix $0<\eta<\min\{\rho_p,1-\rho_p\}$.
The uniform estimate keeps $p_c$ strictly above $p$ for every
$3\le n\le\T_{\lfloor(\rho_p-\eta)d\rfloor}$, and strictly below $p$
at $n=\T_{\lceil(\rho_p+\eta)d\rceil}$. This proves the limit and the
displayed bounds at $p=1/4$ since the intervals widen with $\eta$.
The comparison uses no monotonicity in the side length.
\end{proof}

\section{The first correction}\label{sec:fine}
At height ratio $\rho$, the critical mean balance gives
$p_\rho=(1-\rho)/(2-\rho)$. A site can still reach the threshold when
its mean support falls short by about $\sqrt{\log d}$ standard
deviations. To see the scale on the density axis, observe that the
variance at the mean balance is
$(2-\rho)dp_\rho(1-p_\rho)=dp_\rho$. Increasing $p$ by $\delta$ increases the mean support by
$(2-\rho)d\,\delta$. Put
\[
 \sigma_\rho=\frac{\sqrt{p_\rho}}{2-\rho}
            =\frac{\sqrt{1-\rho}}{(2-\rho)^{3/2}} .
\]

\begin{theorem}[The refined threshold]\label{thm:fine}
Fix $0<\rho_*<1$. Uniformly for $n\ge3$ with $\rho=H(n)/d\le\rho_*$,
\begin{equation}\label{eq:fine}
 p_c(\ZZ_n^d,d)=p_\rho-\sigma_\rho\sqrt{\frac{\log d}{d}}
                 +O_{\rho_*}\bigl((d\log d)^{-1/2}\bigr).
\end{equation}
For every fixed $C>1$ the same estimate holds on boxes of side
$n\ge Cd$. The error constant depends only on $\rho_*$, and the lower
bound on $d$ may also depend on $C$. At $n=d$
the correction is given by Proposition~\ref{prop:poisson}.
\end{theorem}

\begin{figure}[htbp]
\centering
\begin{tikzpicture}
\begin{axis}[width=.61\textwidth,height=4.6cm,xmin=0,xmax=1,ymin=0,ymax=.415,
 xlabel={Height ratio $\rho$},ylabel={$\sigma_\rho$},
 xtick={0,.5,.66666667,1},xticklabels={$0$,$1/2$,$2/3$,$1$},
 ytick={0,.2,.375},yticklabels={$0$,$0.2$,$3/8$},
 tick label style={font=\small}]
\addplot[draw=blue!48!black,very thick,domain=0:1,samples=201,variable=\t]
 ({1-\t^2},{\t/(1+\t^2)^1.5});
\draw[dashed,gray] (axis cs:0,.375) -| (axis cs:.66666667,0);
\addplot[draw=blue!48!black,fill=blue!48!black,only marks,mark=*,mark size=1.8pt]
 coordinates {(.66666667,.375)};
\end{axis}
\end{tikzpicture}
\caption{Coefficient of the correction $-\sigma_\rho\sqrt{\log d/d}$. It
has maximum $2/(3\sqrt3)$ at $\rho=1/2$ and equals $3/8$ at
the height ratio $2/3$. The curve shows
the continuous extension to $\rho=1$.}
\label{fig:refined-coefficient}
\end{figure}
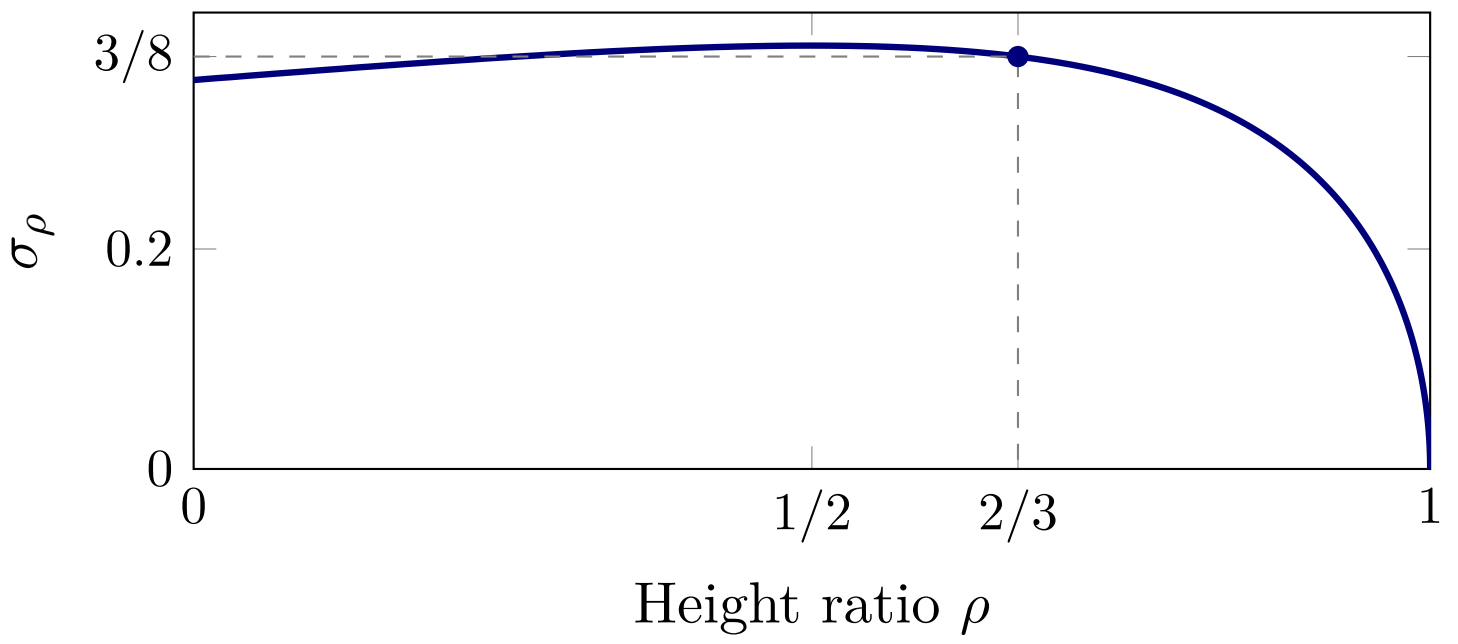

At $\rho=2/3$, corresponding to limiting density $1/4$, the coefficient is $3/8$
(Figure~\ref{fig:refined-coefficient}). Inverting the estimate gives
the height needed at a fixed density.

\begin{corollary}[The height at a prescribed density]\label{cor:fine-height}
For each fixed $0<p<1/2$,
\begin{equation}\label{eq:inverse}
 H(f_p(d))=\frac{1-2p}{1-p}\,d
      -\frac{\sqrt p}{1-p}\sqrt{d\log d}
      +O_p\!\left(\sqrt{\frac d{\log d}}\right).
\end{equation}
In particular $H(f_{1/4}(d))=\tfrac23d-\tfrac23\sqrt{d\log d}
+O(\sqrt{d/\log d})$. At side $n=\T_{\lfloor \rho_p d\rfloor}$ and
density $p$, both graphs fill with probability at least $1-e^{-d}$
for all large $d$.
\end{corollary}

\subsection{The binomial tail}
Let $X\sim\Bin(N,p)$ and $V=Np(1-p)$, where $p$ ranges in a fixed
compact subinterval of $(0,1)$. For any fixed $C$ and
$\sqrt N\le z\le C\sqrt{N\log N}$, uniformly,
\begin{equation}\label{eq:bin}
 \PP(X\ge Np+z)\asymp \frac{\sqrt V}{z}e^{-z^2/(2V)} .
\end{equation}
If $z/\sqrt N\to\infty$, the ratio tends to $1/\sqrt{2\pi}$.
For $z\le\sqrt N$ the tail is bounded below by a positive constant.
We include the argument since the prefactor in~\eqref{eq:bin} will
be used in both estimates (see also~\cite[Chapter~VII]{Feller}).

Put $m=\lceil Np+z\rceil$ and $w=m-Np$. Stirling's formula and a
Taylor expansion of the binomial relative entropy give
$\PP(X=m)=(1+o(1))(2\pi V)^{-1/2}e^{-z^2/(2V)}$. The error in
the exponent is $O(w^3/N^2+(w^2-z^2)/N)=o(1)$. For
$r_j=\PP(X=m+j+1)/\PP(X=m+j)$,
\[
 1-r_j=\frac{w+j+1-p}{(m+j+1)(1-p)} .
\]
The ratios decrease and $1-r_0\sim w/V$, so the tail is at most
$(1+o(1))V/w$ times its first mass. The first $\lfloor V/w\rfloor+1$
terms give a matching bound up to a constant because $1-r_j=O(w/V)$
in that range. When $w/\sqrt N\to\infty$, we have $1-r_j\sim w/V$
uniformly for $0\le j\le T V/w$ for each fixed $T$. Summing those
terms and then letting $T\to\infty$ gives the claimed ratio.
The bound for $z\le\sqrt N$ follows by monotonicity.

\subsection{The lower estimate near the threshold}
Here the gaps between successive neighbour counts have order
$\sqrt{d/\log d}$. Five generations suffice for the rectangle cover.

\begin{proposition}[The lower estimate]\label{prop:fine-lower}
There are absolute constants $C,c>0$ with the following property. Fix
$p_-\in(0,1/2)$ and $A>0$. For all sufficiently large $d$, suppose
$p_-\le p\le1/2$, $H$ is an integer in $[0,d]$ and
$3\le n\le\T_H(d^A)$.
If the deficit in~\eqref{eq:deficit} satisfies
\begin{equation}\label{eq:low-test}
 D_H\ge\sqrt{(2d-H)p(1-p)\log d}+C\sqrt{d/\log d},
\end{equation}
then
\[
 \PP_p(\ZZ_n^d\text{ fills})\le e^{-c(d/\log d)^{5/2}},\qquad
 \frac{\EE_p|A_\infty\setminus A_0|}{n^d}
       =O\bigl((d\log d)^{-1/2}\bigr).
\]
The same upper bounds hold on $[n]^d$.
For fixed $0<p<1/2$, all the local estimates of
Proposition~\ref{prop:first-round} hold throughout \eqref{eq:low-test},
on their respective graphs.
\end{proposition}
\begin{proof}
Take $h=5$ and $g=\lceil8e\sqrt{d/\log d}\rceil$ in
Lemma~\ref{lem:witness}, and write $N=2d-H$.
For sufficiently large $C$, \eqref{eq:low-test} gives
$d-hg\ge H+Np+\sqrt{Np(1-p)\log d}$. A shift of $2h+1$ changes
the quadratic exponent in~\eqref{eq:bin} by $o(1)$. Hence the leaf
bound is $w_0\le(1+o(1))/\sqrt{2\pi d\log d}$ and
$4edw_0\le u=g-2h$. In the first application of
\eqref{eq:recurrence} the bound $(2edw_0/u)^{u/4}\le2^{-u/4}$
suffices. At subsequent levels $\log(2ed/u)=O(\log d)=o(u)$, so
induction gives $w_i\le e^{-c_i g^i}$ for $1\le i\le5$.
Since $g^5/d^2\to\infty$, Proposition~\ref{prop:cover} applies with
$K=c'g^5$ and $t=6$. The empty-background probability is at most
$p+\sum_{i<h}w_i<3/4$. Equations~\eqref{eq:global-fail}
and~\eqref{eq:density-comparison} prove the bounds. They transfer to
the box by deleting periodic edges.

For the local estimates at fixed $p$, use the proof of
Proposition~\ref{prop:first-round} with $G=5g$ and $a=4$. Every site
in $B(\varnothing)\setminus A_0$ lies within $a$ of $E$. For fixed
$M$ we have $M\le G=o(d)$, and the exceptional probability $e^{-K/3}$
is smaller than $e^{-Cd}$ for every fixed $C$. The same component and
two-neighbour calculations therefore apply with these parameters.
\end{proof}

Any fixed $h\ge5$ gives an exponent $c_h(d/\log d)^{h/2}$ in the
first bound, on increasing $C$ in~\eqref{eq:low-test}. Only the test
radius and the separation $2(h+1)$ in Lemma~\ref{lem:slabs} change.

\subsection{Growth from outward neighbours}
\begin{figure}[htbp]
\centering
\begin{tikzpicture}[x=6mm,y=6mm,>=latex,font=\footnotesize]
  \begin{scope}
    \coordinate (x) at (0,0);
    \coordinate (a) at (2,1.25);
    \coordinate (b) at (2,-1.25);
    \coordinate (u) at (4,2);
    \coordinate (z) at (4,0);
    \coordinate (v) at (4,-2);
    \foreach \a/\b in {x/a,x/b,a/u,b/v}
      \draw[->,black!55,shorten <=2.8pt,shorten >=2.8pt] (\a)--(\b);
    \foreach \a in {a,b}
      \draw[->,blue!60!black,thick,shorten <=2.8pt,shorten >=2.8pt] (\a)--(z);
    \foreach \p in {x,a,b,u,v}
      \fill[black!70] (\p) circle (2.15pt);
    \fill[blue!70!black] (z) circle (2.45pt);
    \node[left] at (x) {$0$};
    \node[above left] at (a) {$e_1$};
    \node[below left] at (b) {$e_2$};
    \node[right] at (u) {$2e_1$};
    \node[right] at (z) {$e_1+e_2$};
    \node[right] at (v) {$2e_2$};
    \node[below] at (2.2,-2.5) {One seed shared by two tests};
  \end{scope}
  \begin{scope}[xshift=7.0cm]
    \coordinate (w) at (0,0);
    \foreach \h in {1.95,.65,-.65,-1.95}
      \draw[black!50] (w)--(2.2,\h);
    \fill[black!70] (w) circle (2.15pt);
    \node[left] at (w) {$w$};
    \foreach \h in {1.95,.65}
      \fill[blue!70!black] (2.2,\h) circle (2.7pt);
    \foreach \h in {-.65,-1.95}
      \draw[orange!75!black,line width=1.1pt,fill=white] (2.2,\h) circle (2.7pt);
    \node[anchor=west,align=left] at (2.65,1.3) {initially\\infected};
    \node[anchor=west,align=left] at (2.65,-1.3) {infected in\\a later round};
    \node[below] at (2.2,-2.5) {Two disjoint sources of support};
  \end{scope}
\end{tikzpicture}
\caption{Both children of the root use the seed at $e_1+e_2$. Along
outward paths a grandchild contributes to the tests of at most
$2\kappa$ children, by~\eqref{eq:backward}. On the right, a site with
threshold four has two initially infected children. Two more children
that were initially healthy bring it to its threshold.}
\label{fig:outward}
\end{figure}
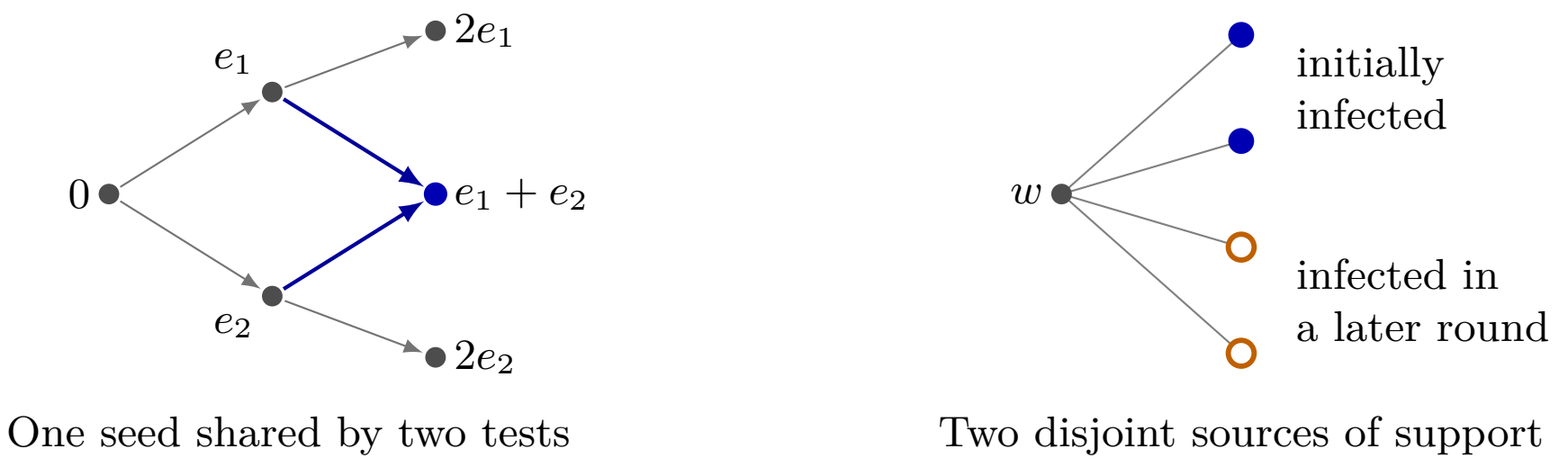

We turn to the upper bound. With the notation of
Section~\ref{sec:box}, the variance of the initial support at a site
of $\mathcal M(k,\ell)$, including a corner, is
\begin{equation}\label{eq:variance}
 (2k+\ell)pq=ps+p\Delta .
\end{equation}
Thus near the leading threshold it is $ps+o(s)$. An initially healthy
neighbour that falls just short in the first round can be brought into
the process by its own newly infected neighbours. We count these
infections going outwards from a fixed root. Two children may see the
same seed, as in Figure~\ref{fig:outward}, but this overlap is bounded.

\begin{lemma}[Local filling]\label{lem:fine-local}
Fix $p_-\in(0,1/2)$ and $\alpha>0$. Let $p\in[p_-,1/2]$, $s=k+\ell$,
$k\ge\alpha s$, $N=2k+\ell$ and $V=Np(1-p)$. Consider a torus of side
at least three, or a box of side lengths at least two with exterior
seeds, in each case with threshold $k$ and $\ell$ private bits per
site. There are constants
$C_b,c>0$, depending only on $p_-,\alpha$, such that if $s$ is large
and $k\le Np+\sqrt{V(\log s-C_b)}$, then for some
$t=O(\log\log\log s)$ with $t\le H(s)+O(1)$,
\begin{equation}\label{eq:weak-local}
 \sup_x\PP_p(x\text{ is healthy at time }t)\le e^{-c(\log s)^3}.
\end{equation}
In particular the hypothesis holds if $\Delta\le\sqrt{ps(\log s-C_b)}$.
\end{lemma}

\begin{proof}
Write $L=\log s$, $a=\sqrt{V/L}$ and
$r_0=\lfloor Np+\sqrt{V(L-C_b)}\rfloor\ge k$. The last assertion
of the lemma follows from $\Delta=k-Np$ and the identity $V=ps+p\Delta$
in~\eqref{eq:variance}.

A child lies on the next sphere about a fixed root $x$, using ambient
distance on a box. At a site $v$ let $Z_b(v)$ count its initially
healthy updating children $w$ for which the private bits and the
initially infected children of $w$ together supply at least $r_0-b$.
Exterior seeds count toward initial support. For integers
$0=b_0\le b_1\le\cdots\le b_T$, impose the tests
\begin{equation}\label{eq:sphere-test}
 Z_{b_{i-1}}(v)\ge b_i\qquad
    (1\le i\le T,\ \ v\in S(x,T-i)\text{ an updating site}).
\end{equation}
If the tests pass, each tested vertex has at least $b_i$ initially
healthy children infected by round $i$. Indeed, a counted child has
$r_0-b_{i-1}$ initial units of support. The preceding test supplies
$b_{i-1}$ further children, all initially healthy, so it reaches
threshold $k$ (Figure~\ref{fig:outward}). For $i=1$ no further children
are needed. Each child's health and private bits enter only its own
indicator. Any grandchild bit enters at most $2\kappa$ child tests,
by~\eqref{eq:backward} applied from that grandchild. Thus if the target
in each test is at most a fixed fraction $1-\eta$ of its mean and
every mean is at least $c_0a$, \eqref{eq:read} bounds the total
probability of a failed test by $(2k+1)^T e^{-c_1a}$.

For $T=O(\log\log\log s)$, a child's outward and private support has
law $\Bin(N-j,p)$ where $0\le j\le\kappa T$. Write $V_j=(N-j)pq$.
At $b=\lfloor a\theta\rfloor$ the standardized excess satisfies
\[
 \frac{r_0-b-(N-j)p}{\sqrt{V_j}}
 =\sqrt{L-C_b}-\frac{\theta}{\sqrt L}
     +O\!\left(\frac{T+1}{\sqrt s}+\frac{T\sqrt L}{s}\right).
\]
Uniformly for $0\le\theta\le 6\log L$, the right side tends to infinity
and its square equals $L-C_b-2\theta+o(1)$. Equation~\eqref{eq:bin}
therefore gives a tail of
$(1+o(1))e^{C_b/2+\theta}/\sqrt{2\pi sL}$.
Even at a corner, $v$ has at least $k-\kappa T$ updating children.
Each starts healthy with probability $q$, independently of its support.
Since $k\ge\alpha s$ and $a\asymp\sqrt{s/L}$, we can choose $C_b$ so that
\begin{equation}\label{eq:gains}
 \EE Z_{\lfloor a\theta\rfloor}(v)\ge4a e^\theta .
\end{equation}

Start at $\theta_0=0$ and take
$\theta_{j+1}=\min\{3\log L,\,2e^{\theta_j}\}$ until the cap is reached
at $J$. The uncapped positive terms at least double, so
$J=O(\log\log\log s)$. Since $2e^z\ge 2^z$ for $z\ge0$, they also
dominate the successive towers of twos and $J\le H(s)+O(1)$. Set $T=J+1$, $b_i=\lfloor a\theta_i\rfloor$
for $i\le J$ and $b_T=\lfloor a L^2\rfloor$. By~\eqref{eq:gains}
each target is at most half its mean, including the last one, whose
mean is at least $4aL^3$. Finally $r_0-Np\le aL$. Hoeffding's
inequality gives probability at most $e^{-c_2L^3}$ that the root's
initial support is less than $r_0-b_T$. Otherwise the $b_T$ initially
healthy children supplied by \eqref{eq:sphere-test} infect the root
by $t=T+1$. Since $T\log(2k+1)+L^3=o(a)$, this proves
\eqref{eq:weak-local}.
\end{proof}

On the $k$-dimensional torus at threshold $k$, suppose the one-site
survival probability at time $t$ is at most $e^{-Q}$ everywhere. For an integer $1\le u\le k/4$ with
$Q\ge2(1+u\log(2k))$,
\[
 \PP_p(x\text{ is healthy at }t+u)\le
 \exp\!\left\{-\frac{Q(k/4)^u}{2u!(2k+1)^t}\right\}.
\]
Indeed a surviving site needs at least $k$ healthy neighbours one
round earlier. On sphere $i$ it has at least $k-2i$ outward choices
and each site on the next sphere has at most $2(i+1)$ parents.
Iteration gives at least $P=(k/4)^u/u!$ survivors on sphere $u$ at
time $t$. Their mean is at most $(2k)^u e^{-Q}$ and each input bit
enters at most $(2k+1)^t$ indicators. Now \eqref{eq:read} applies
since $e(2k)^u e^{-Q}/P\le e^{-Q/2}$. With $u=t+2$ and $Q=c(\log s)^3$
from Lemma~\ref{lem:fine-local} this yields
$\PP_p(x\text{ is healthy at }2t+2)\le e^{-s^{7/4}}$, as the exponent
in the display is at least $s^{2-o(1)}$.

\subsection{A sufficient height}
\begin{proposition}[The upper estimate]\label{prop:fine-up}
Fix $p_-\in(0,1/2)$. There are $C_b,C_s>0$ such that, uniformly in
$p\in[p_-,1/2]$ and for all sufficiently large $d$, the box and the
torus fill with probability at least $1-e^{-d}$ whenever
\begin{equation}\label{eq:high-test}
 n\ge\T_{\widehat H_p},\qquad
 \widehat H_p=\max\{0,\lceil A_p\rceil\}+H(d)+2,\qquad
 A_p=\frac{(1-2p)d-\sqrt{pd(\log d-C_b)}}{1-p}.
\end{equation}
If $A_p\le0$, this probability bound holds for every torus of side
$n\ge3$ and every box of side $n\ge C_s d$.
\end{proposition}
\begin{proof}
Choose $C_b$ in Lemma~\ref{lem:fine-local} for $\alpha=p_-^2/100$. Take
$W(s)=\sqrt{ps(\log s-C_b)}$ for $s\ge s_0$ and extend it linearly to
zero on $[0,s_0]$. For sufficiently large $s_0$, $W(s)/s$ is
nonincreasing. At every site of a cell with exterior seeds, the
bound~\eqref{eq:weak-local} verifies~\eqref{eq:start-local}, since
$(\log s)^3\gg t\log(2k+1)$. Now $J=\max\{0,\lceil A_p\rceil\}$
in~\eqref{eq:potential}. Lemma~\ref{lem:base} and
Proposition~\ref{prop:growth} give $R_J+2\le\T_{\widehat H_p}$ and
prove~\eqref{eq:high-test}. When $A_p\le0$, all cells start at a base
pair and $R_0+2\le C_s d$.

Adding periodic edges transfers the box bound to the torus. If
$A_p\le0$, Lemma~\ref{lem:fine-local} also applies with $k=s=d$ on every
torus. Equation~\eqref{eq:torus-path} gives failure probability at most $e^{-d}$
uniformly in its side length.
\end{proof}

\subsection{Reading the threshold from the two estimates}
Set $L=\log d$, $\omega=\sqrt{d/L}$, $H=H(n)$ and $S=\sqrt{pdL}$.
For $p$ in a fixed compact subinterval of $(0,1/2]$ both comparisons
are governed by
\begin{equation}\label{eq:balance}
 F(H,p)=(1-2p)d-(1-p)H-S=D_H-S.
\end{equation}
For a sufficiently large fixed $K$, $F\ge K\omega$ implies
\eqref{eq:low-test}. Indeed $(2d-H)p(1-p)=p(d+D_H)$, and
$D-\sqrt{p(d+D)L}$ is increasing for $D\ge0$ when $d$ is large.
At $D=S+K\omega$ its value is $K\omega-O(L)$, where $L=o(\omega)$.
Conversely, $F\le-K\omega$ gives $A_p\le H-2\omega$ because
$(1-p)(A_p-H)=D_H-\sqrt{pd(L-C_b)}$ and
$S-\sqrt{pd(L-C_b)}=O(\omega)$. If $A_p>0$, then
$\widehat H_p\le H-1$ since $H(d)+4=o(\omega)$, and
$n>\T_{H-1}$ gives~\eqref{eq:high-test}.
If $A_p\le0$, Proposition~\ref{prop:fine-up} applies on every torus
and on boxes of side at least $C_s d$.

\begin{proof}[Proof of Theorem~\ref{thm:fine}]
Let $\rho=H/d\le\rho_*<1$. At
$p=p_\rho-\sigma_\rho\sqrt{L/d}$ we have
$D_H=(2-\rho)d(p_\rho-p)=\sqrt{p_\rho dL}$ and $F=O_{\rho_*}(L)$.
Changing $p$ by $\pm\Lambda/\sqrt{dL}$ changes $F$ by
$\mp(2-\rho+o(1))\Lambda\omega$, uniformly. The two test densities
$p^-,p^+$ stay in a compact subinterval of $(0,1/2]$. Taking $\Lambda$
large in the preceding comparisons proves the result for all tori and
for boxes with $n\ge C_s d$.

To extend the box result to every fixed $C>1$, it remains to consider
$Cd\le n\le C_s d$. Here $H=o(\omega)$, so $A_{p^+}\le0$.
With $W$ from Proposition~\ref{prop:fine-up}, this gives
$W(k)/k\ge W(d)/d\ge1-2p^+$ for $k\le d$. Hence for all sufficiently
large $k$ the $k$-cells meet the hypothesis of Lemma~\ref{lem:fine-local}.
Put $m=n-2$. Since $k\le m$, $k\log m\le m\log k$ for $k\ge3$,
and $Q/t\gg\log k$ in \eqref{eq:weak-local}, Lemma~\ref{lem:path}
bounds the probability that a given such cell fails by $e^{-2m}$.
Their total contribution is at most $3^d e^{-2m}=o(1)$.
For each of the finitely many smaller dimensions use
Lemma~\ref{lem:fixed-face} at a fixed $p_0<1/2$ with
$-C\log(1-p_0)>\log2$. Since $p^+\to1/2$, the sum over these cells
also tends to zero, as in the proof of~\eqref{eq:linear-box}.
The bound at $p^-$ transfers from the torus by edge deletion.
\end{proof}

\begin{proof}[Proof of Corollary~\ref{cor:fine-height}]
At fixed $p$, equation $F(H,p)=0$ has the exact solution
$H_0=\rho_p d-\sqrt{pdL}/(1-p)$. For a large fixed $K'$,
$F\ge K\omega$ at every height below $H_0-K'\omega$.
Thus all sides up to $\T_{\lfloor H_0-K'\omega\rfloor}$ have filling
probability below $1/2$. At $n=\T_{\lceil H_0+K'\omega\rceil}$ the
opposite comparison gives filling probability at least $1-e^{-d}$.
Continuity in $p$ on finite graphs now brackets the first crossing and
proves~\eqref{eq:inverse}. At $H=\lfloor\rho_pd\rfloor$ we have
$F=-\sqrt{pdL}+O(1)$, giving the same filling bound on both graphs.
\end{proof}

\section*{Acknowledgements}
The authors thank Robert Morris for posing the torus question at the
Tsinghua Probability Seminar in October 2026. The research of Yuval Peres
was supported by National Natural Science Foundation of China grant
RFIS-W2531011.
During the preparation of this work, the authors used several large language models. The authors reviewed and verified all AI-assisted content and take full responsibility for the final manuscript.


\begin{thebibliography}{99}
\setlength{\itemsep}{2pt plus 1pt minus 1pt}

\bibitem{AL}
M.~Aizenman and J.~L.~Lebowitz,
\href{https://doi.org/10.1088/0305-4470/21/19/017}{\emph{Metastability effects in bootstrap percolation}},
J. Phys. A: Math. Gen. \textbf{21} (1988), no.~19, 3801--3813.

\bibitem{BBDCM}
J.~Balogh, B.~Bollob\'as, H.~Duminil-Copin and R.~Morris,
\href{https://doi.org/10.1090/S0002-9947-2011-05552-2}{\emph{The sharp threshold for bootstrap percolation in all dimensions}},
Trans. Amer. Math. Soc. \textbf{364} (2012), no.~5, 2667--2701.
\href{https://arxiv.org/abs/1010.3326}{arXiv:1010.3326}.

\bibitem{BBM3}
J.~Balogh, B.~Bollob\'as and R.~Morris,
\href{https://doi.org/10.1214/08-AOP433}{\emph{Bootstrap percolation in three dimensions}},
Ann. Probab. \textbf{37} (2009), no.~4, 1329--1380.

\bibitem{BBM}
J.~Balogh, B.~Bollob\'as and R.~Morris,
\href{https://doi.org/10.1017/S0963548308009322}{Combin. Probab. Comput.
\textbf{18} (2009), nos.~1--2, 17--51}.
\href{https://arxiv.org/pdf/math/0702373\#page=37}{arXiv:math/0702373, Problem~1}.

\bibitem{Book}
A.~Bou-Rabee, L.~Levine and Y.~Peres,
\emph{Discrete Laplacian Growth: Sandpiles, Rotors, and DLA},
book manuscript, 10 October 2026, \S2.11, Lemma~2.45.

\bibitem{CC}
R.~Cerf and E.~N.~M.~Cirillo,
\href{https://doi.org/10.1214/aop/1022874817}{\emph{Finite size scaling in three-dimensional bootstrap percolation}},
Ann. Probab. \textbf{27} (1999), no.~4, 1837--1850.

\bibitem{CM}
R.~Cerf and F.~Manzo,
\href{https://doi.org/10.1016/S0304-4149(02)00124-2}{\emph{The threshold regime of finite volume bootstrap percolation}},
Stochastic Process. Appl. \textbf{101} (2002), no.~1, 69--82.

\bibitem{CLR}
J.~Chalupa, P.~L.~Leath and G.~R.~Reich,
\href{https://doi.org/10.1088/0022-3719/12/1/008}{\emph{Bootstrap percolation on a Bethe lattice}},
J. Phys. C: Solid State Phys. \textbf{12} (1979), no.~1, L31--L35.

\bibitem{CEGK2}
M.~Collares, J.~Erde, A.~Geisler and M.~Kang,
\href{https://doi.org/10.1016/j.ejc.2026.104347}{\emph{Majority bootstrap percolation on the permutahedron and other
high-dimensional graphs}},
European J. Combin. \textbf{134} (2026), article 104347.
\href{https://arxiv.org/abs/2507.06597}{arXiv:2507.06597}.

\bibitem{Feller}
W.~Feller,
\emph{An Introduction to Probability Theory and Its Applications}, Vol.~I,
3rd ed., John Wiley \& Sons, New York, 1968.

\bibitem{Finner}
H.~Finner,
\href{https://doi.org/10.1214/aop/1176989534}{\emph{A generalization of H\"older's inequality and some probability inequalities}},
Ann. Probab. \textbf{20} (1992), no.~4, 1893--1901.

\bibitem{GLSS}
D.~Gavinsky, S.~Lovett, M.~Saks and S.~Srinivasan,
\href{https://doi.org/10.1002/rsa.20532}{\emph{A tail bound for read-$k$ families of functions}},
Random Structures Algorithms \textbf{47} (2015), no.~1, 99--108.

\bibitem{Hoeffding}
W.~Hoeffding,
\href{https://doi.org/10.1080/01621459.1963.10500830}{\emph{Probability inequalities for sums of bounded random variables}},
J. Amer. Statist. Assoc. \textbf{58} (1963), no.~301, 13--30.

\bibitem{Holroyd}
A.~E.~Holroyd,
\href{https://doi.org/10.1007/s00440-002-0239-x}{\emph{Sharp metastability threshold for two-dimensional bootstrap
percolation}},
Probab. Theory Related Fields \textbf{125} (2003), no.~2, 195--224.

\bibitem{MorrisSurvey}
R.~Morris,
\href{https://doi.org/10.1017/9781108332699.008}{\emph{Monotone cellular automata}},
in \emph{Surveys in Combinatorics 2017}
(A.~Claesson, M.~Dukes, S.~Kitaev, D.~Manlove and K.~Meeks, eds.),
London Math. Soc. Lecture Note Ser., vol.~440, Cambridge Univ. Press,
2017, pp.~312--371.
\href{https://w3.impa.br/~rob/BCC_survey.pdf}{Author's version}.

\bibitem{Schonmann}
R.~H.~Schonmann,
\href{https://doi.org/10.1214/aop/1176989923}{\emph{On the behavior of some cellular automata related to bootstrap
percolation}},
Ann. Probab. \textbf{20} (1992), no.~1, 174--193.

\end{thebibliography}
\end{document}